\documentclass[11pt,a4paper]{article}

\usepackage{a4wide}
\usepackage{graphicx}
\usepackage{epstopdf}
\usepackage{epsfig}
\usepackage[usenames,dvipsnames]{color}
\usepackage{courier}
\usepackage{capt-of}
\usepackage{hyperref}
\usepackage{underscore}	% for the underscore
\usepackage{amsmath}
\usepackage{amsthm}
\usepackage{verbatim}
\usepackage{caption}
\usepackage{listings}
\usepackage{pstricks}
\usepackage{booktabs}
\usepackage{multirow}
\usepackage{bm}
\usepackage{amssymb}
\usepackage[T1]{fontenc}
\usepackage[utf8]{inputenc}
\usepackage{setspace}
\usepackage[english]{babel}
\usepackage{grffile}
\usepackage{algorithm}
\usepackage{algorithmicx}
\usepackage{ifthen}
\usepackage{algpseudocode}
\usepackage{subfig}
\usepackage{subfloat}
\usepackage{amsfonts}
\usepackage{fancyhdr}
\usepackage{varwidth}
\usepackage{appendix}
\usepackage{tikz}

\makeatletter
\def\cleardoublepage{\clearpage\if@twoside \ifodd\c@page\else
    \hbox{}
    \thispagestyle{plain}
    \newpage
    \if@twocolumn\hbox{}\newpage\fi\fi\fi}
\makeatother \clearpage{\pagestyle{plain}\cleardoublepage}

\newcommand{\EndProofMarker}{$\Box$}

\usepackage[a4paper, margin=2.3cm]{geometry}

\newtheorem{remark}{Remark}[section]

\newcommand{\param}{{\bm{\mu}}}
\newcommand{\paramspace}{\mathcal{D}}
\newcommand{\Nh}{N_h}
\newcommand{\real}{\mathbb{R}}
\newcommand{\realh}{\real^{\Nh}}
\newcommand{\realhh}{\real^{\Nh\times \Nh}}

\newcommand{\realN}{\real^{N}}
\newcommand{\realNN}{\real^{N\times N}}

\newcommand{\uu}{u}

\newcommand{\umu}{u(\param)}

\newcommand{\nonlinear}{\mathcal{N}}
\newcommand{\nonlinearMatrix}{\mathbf{N}}

\newcommand{\vecuh}{\mathbf{\uu}}
\newcommand{\vecuhi}[1]{\vecuh_{#1}}
\newcommand{\vecuhimu}[1]{\vecuhi{#1}^\param}
\newcommand{\vecfh}{\mathbf{f}}
\newcommand{\vecfhmu}{\vecfh(\param)}

\newcommand{\rbbasis}{\xi}
\newcommand{\rbvecbasis}{\bm{\rbbasis}}
\newcommand{\rbV}{\mathbf{V}}

\newcommand{\epod}{\varepsilon_{\operatorname{POD}}}

\newcommand{\vecf}{\mathbf{f}}
\newcommand{\vecfmu}{\vecf(\param)}
\newcommand{\vecfq}{\vecf^q}
\newcommand{\nonlinearMatrixmu}{\nonlinearMatrix[\vecuh, \param]}
\newcommand{\nonlinearMatrixuNmu}{\nonlinearMatrix[\rbV\rbvecu, \param]}
\newcommand{\nonlinearMatrixq}{\nonlinearMatrix^q}

\newcommand{\rbvecu}{\mathbf{\uu}_N}
\newcommand{\rbvecumu}{\rbvecu(\param)}
\newcommand{\rbNonlinearMatrix}{\mathbf{N}_N}
\newcommand{\rbNonlinearMatrixmu}{\rbNonlinearMatrix[\rbvecu, \param]}
\newcommand{\rbvecf}{\mathbf{f}_N}
\newcommand{\rbvecfq}{\rbvecf^q}
\newcommand{\rbvecfmu}{\rbvecf(\param)}
\newcommand{\rbNonlinearMatrixq}{\rbNonlinearMatrix^q}
\newcommand{\Qn}{Q_n}
\newcommand{\thetan}{\theta_n^q}
\newcommand{\thetann}[1]{\theta_n^{#1}}
\newcommand{\thetanmu}{\thetan(\rbvecu;\param)}
\newcommand{\Qf}{Q_f}
\newcommand{\thetaf}{\theta_f^q}
\newcommand{\thetaff}[1]{\theta_f^{#1}}
\newcommand{\thetafmu}{\thetaf(\param)}

\newcommand{\ns}{n_s}

\newcommand{\snap}{\mathbf{S}}

\newcommand{\snaps}{\mathbf{u}}
\newcommand{\snapi}{\snaps(\param_i)}
\newcommand{\setsnapi}{ \{\snapi\}_{i=1}^{\ns} }
\newcommand{\svdU}{\mathbf{U}}
\newcommand{\svdSigma}{\bm{\Sigma}}
\newcommand{\svdZ}{\mathbf{Z}}
\newcommand{\realnsns}{\real^{\ns\times\ns}}
\newcommand{\realhns}{\real^{\Nh\times\ns}}

\newcommand{\wl}{W^{(l)}}
\newcommand{\wL}{W^{(L)}}
\newcommand{\bl}{\mathbf b^{(l)}}
\newcommand{\bL}{\mathbf b^{(L)}}

\newcommand{\xxl}{\mathbf x^{(l)}}
\newcommand{\xxL}{\mathbf x^{(L)}}
\newcommand{\XX}{\mathbf X}
\newcommand{\XXl}{\mathbf X^{(l)}}
\newcommand{\yyl}{\mathbf y^{(l)}}
\newcommand{\yyL}{\mathbf y^{(L)}}
\newcommand{\YY}{\mathbf Y}
\newcommand{\YYl}{\mathbf Y^{(l)}}
\newcommand{\yylm}{\mathbf y^{(l-1)}}
\newcommand{\sigmal}{\sigma^{(l)}}
\newcommand{\sigmaRB}{\sigma_{\operatorname{RB}}}
\newcommand{\relu}{\text{ReLU}}
\newcommand{\sigm}{\text{sigmoid}}
\newcommand{\Nin}{N_{\text{in}}}
\newcommand{\Nout}{N_{\text{out}}}
\newcommand{\Ninl}{\Nin^{(l)}}
\newcommand{\Noutl}{\Nout^{(l)}}
\newcommand{\realNin}{\real^{\Nin}}
\newcommand{\realNout}{\real^{\Nout}}
\newcommand{\realNinl}{\real^{\Ninl}}
\newcommand{\realNoutl}{\real^{\Noutl}}
\newcommand{\realNoutlNinl}{\real^{\Noutl \times \Ninl}}

\newcommand{\loss}{\mathcal L}
\newcommand{\thetal}{\theta^{(l)}}
\newcommand{\MSE}{\text{MSE}}

\newcommand{\Output}[1]{Output \##1}
\newcommand{\Input}[1]{Input \##1}

\newcommand{\Dmu}[1]{\mathcal{O}_{#1}^\param}

\newcommand{\Psetin}{\mathcal P_\text{in}}
\newcommand{\Psetout}{\mathcal Q_\text{out}}
\newcommand{\pp}{\vec{p}}
\newcommand{\ppo}{\vec{q}}
\newcommand{\Ns}{N_\text{s}}
\newcommand{\xxi}{\mathbf x_i}
\newcommand{\yyi}{\mathbf y_i}

\newcommand{\latent}{\boldsymbol \xi}

\newcommand{\vecb}{\vec{b}}

\newcommand{\vecx}{\vec{x}}

\newcommand{\diffCoeffmu}{\alpha( \param )}
\newcommand{\feAhNa}{\mathbf{A}}
\newcommand{\feAhNamu}{\feAhNa( \param )}
\newcommand{\feAhNaq}{\feAhNa_q}
\newcommand{\thetaAhNamu}[1]{\theta^a_{#1}( \param )}
\newcommand{\QaNa}{Q_a}

\newcommand{\feFNa}{\mathbf{F}}
\newcommand{\feFNamu}{\feFNa( \param )}
\newcommand{\QfNa}{Q_f}
\newcommand{\feFNaq}{\feFNa_q}
\newcommand{\thetaFNamu}[1]{\theta^f_{#1}( \param )}

\newcommand{\feTNa}{\mathbf{T}}
\newcommand{\feTNamu}{\feTNa( \param )}

\newcommand{\Nhu}{N^u_{h}}
\newcommand{\Nhp}{N^p_{h}}
\newcommand{\nsvecu}{\mathbf{u}}
\newcommand{\nsvecp}{\mathbf{p}}
\newcommand{\nsvecfu}{\mathbf{f}^u}
\newcommand{\nsvecfumu}{\nsvecfu(\param)}

\newcommand{\nsu}{\vec{u}}
\newcommand{\nsumu}{\nsu(\param)}
\newcommand{\nsp}{p}
\newcommand{\nspmu}{\nsp(\param)}

\newcommand{\rey}{Re}
\newcommand{\nsnu}{\nu}

\newcommand{\nsD}{\mathbf{D}}
\newcommand{\nsDmu}{\nsD(\param)}
\newcommand{\nsK}{\mathbf{K}}
\newcommand{\nsKmu}{\nsK(\param)}
\newcommand{\nsC}{\mathbf{C}}
\newcommand{\nsCmu}{\nsC(\nsvecu;\param)}

\newcommand{\feBh}{\mathbf{B}}

\newcommand{\nsrbvecfu}{\mathbf{f}^u_N}
\newcommand{\nsrbvecfumu}{\nsrbvecfu(\param)}
\newcommand{\nsrbD}{\mathbf{D}_N}
\newcommand{\nsrbDmu}{\nsrbD(\param)}
\newcommand{\nsrbK}{\nsK_N}
\newcommand{\nsrbKmu}{\nsrbK(\param)}

\newcommand{\rbVu}{\mathbf{V}_{u}}

\newcommand{\Nu}{N_u}

\newcommand{\nsrbvecun}{\mathbf{u}_N}
\newcommand{\nsrbvecunmu}{\nsrbvecun(\param)}
\newcommand{\nsCunmu}{\mathbf{C}_N(\rbVu\nsrbvecun; \param)}
\newcommand{\nsCj}{\mathbf{C}_j}

\newcommand{\QkNs}{Q_k}
\newcommand{\thetaKhNsmu}[1]{\theta^k_{#1}( \param )}
\newcommand{\feKhNsq}{\nsK_q}

\newcommand{\QfuNs}{Q_f^u}
\newcommand{\thetaFhNsmu}[1]{\theta^{f^u}_{#1}( \param )}
\newcommand{\feFhNsq}{\nsvecfu_q}
\newcommand{\unj}{u_N^j}
\newcommand{\uni}[1]{u_N^{#1}}
\newcommand{\QnNs}{Q_n}

\newcommand{\archNaOne}{\mathcal{A}^1_{\operatorname{na}}}
\newcommand{\archNaTwo}{\mathcal{A}^2_{\operatorname{na}}}
\newcommand{\archNaThree}{\mathcal{A}^3_{\operatorname{na}}}

\newcommand{\archNsOne}{\mathcal{A}^1_{\operatorname{ns}}}
\newcommand{\archNsTwo}{\mathcal{A}^2_{\operatorname{ns}}}
\newcommand{\archNsThree}{\mathcal{A}^3_{\operatorname{ns}}}
\newcommand{\archNsFour}{\mathcal{A}^4_{\operatorname{ns}}}

\begin{document}

\title{Data driven approximation of parametrized PDEs by Reduced Basis and Neural Networks}
%\title{Approximation of parameterized PDEs by Concurrent Reduced Basis Method and Neural Networks}
\author{Niccol\`o Dal Santo$^1$, Simone Deparis$^1$, Luca Pegolotti$^1$}

\maketitle
\begin{center}
	{$^1$ SCI-SB-SD, \'{E}cole Polytechnique F\'{e}d\'{e}rale de Lausanne (EPFL), \\
		Station 8, 1015 Lausanne, Switzerland. \\
	}
\end{center}

\date{}

\begin{abstract}
We are interested in the approximation of partial differential equations with a data-driven approach based on the reduced basis method and machine learning.
We suppose that the phenomenon of interest can be modeled by a parametrized partial differential equation, but that the value of the physical parameters is unknown or difficult to be directly measured.
Our method allows to estimate fields of interest, for instance temperature of a sample of material or velocity of a fluid, given data at a handful of points in the domain. We propose to accomplish this task with a neural network embedding a reduced basis solver as exotic activation function in the last layer. The reduced basis solver accounts for the underlying physical phenomenonon and it is constructed from snapshots obtained from randomly selected values of the physical parameters during an expensive offline phase. The same full order solutions are then employed for the training of the neural network. As a matter of fact, the chosen architecture resembles an asymmetric autoencoder in which the decoder is the reduced basis solver and as such it does not contain trainable parameters.
The resulting latent space of our autoencoder includes parameter-dependent quantities feeding the reduced basis solver, which -- depending on the considered partial differential equation -- are the values of the physical parameters themselves or the affine decomposition coefficients of the differential operators.

\end{abstract}

\section{Introduction}
During the last decade, machine learning (ML) techniques have gained considerable attention for their
ability to exploit data abundance in problems that are typically difficult to be solved through
classic algorithmic paradigms. This success is partly motivated by the recent exponential
growth in the size of available data, which is due to the technological advancements in consumer
devices, and by the increase in computational power of parallel architectures, such as
supercomputers and graphics processing units (GPUs). Deep learning (DL) \cite{goodfellow2016deep,bengio2009learning}
currently represents one of the most active and promising areas of research in machine learning. Notable applications of DL
include classification tasks, for instance image recognition \cite{krizhevsky2012imagenet,mordvintsev2015inceptionism,wang2017residual}, text categorization
\cite{lai2015recurrent,kim2014convolutional} and natural language processing \cite{young2018recent,collobert2011natural}, and regression tasks \cite{specht1991general}, which play an important role e.g. in computer vision \cite{lathuiliere2018comprehensive}, object detection \cite{szegedy2013deep},
stock prediction \cite{ding2015deep}, and cancer detection \cite{hu2018deep}.

In the context of the numerical solution of partial differential equations (PDEs), DL
has only lately become an alternative to standard approximation techniques. This is arguably due to the reluctancy
of the community to adopt these algorithms because of the lack of a solid theoretical foundation, which instead
characterizes classical numerical methods -- e.g. the finite element (FE) method -- and provides
error bounds and stability conditions. In fact, the approximation properties of neural networks (NNs), which are the basis of every DL algorithm,
have been well-known for at least thirty years; for example, in \cite{university1988continuous, cybenko1989approximation,mhaskar1993approximation}
the authors show that feed-forward neural networks with one hidden layer are universal approximators. However, these
results do not address fundamental aspects regarding the application of neural networks in practice;
for example, they do not provide guidelines on how to design optimal architectures retaining a sufficient
level of expressiveness with a minimal number of parameters and they do not focus on the
\textit{trainability} of such algorithms.
Some of the recent efforts oriented towards defining a mathematical theory of
neural networks and which give insights on the aformentioned critcal issues are e.g. \cite{han2018mean,bolcskei2019optimal,petersen2018topological,hanin2017universal,hanin2018neural,hanin2018start}.

In this paper, we address the problem of approximating the solution of parametrized
PDEs by combining the reduced basis (RB) method \cite{quarteroni2015reduced,hesthaven2015certified} and NNs.
We refer the reader to \cite{lee2018model,hesthaven2018non,wang2019non} for other
applications of DL to reduced order modelling with the RB method.
We set ourselves in the scenario in which we are given data at some input physical points -- e.g. the value of the
variable of interest -- and we are interested in determining the solution itself (or functions thereof) at output locations or
even in the whole domain.
Moreover, we suppose that the values of the physical parameters
characterizing the PDE be either unknown or difficult to be measured (e.g. with
non-invasive procedures). The key idea of our approach is to design a NN to learn
the input-output mapping by exploiting the knowledge of the underlying physical phenomenon. This concept
has already been introduced in the works by E et al. \cite{E2018}, Raissi et al. \cite{raissiI2017physics,raissiII2017physics}, and
Schwab and Zech \cite{schwab2019}.
In these cases, the authors propose to include the knowledge of the PDE in the
functional that is minimized during training: in the former work, the
loss function to be minimized is identified as the energy functional of the
corresponding differential problem, whereas in the latter the classic mean squared
error loss function is modified to account for the residual of the PDE (in strong form)
and the initial and boundary conditions.
The novelty of this paper consists in including the physical knowledge directly in the NN
by considering as output layer a RB solver that, given
the input of the previous layers, assembles and solves the reduced system, projects
the reduced solution onto the original full-order space, and computes the solution at the desired
output locations. We refer to this type of architecture as PDE-aware deep neural network (PDE-DNN).
This can be interpreted as a standard multi-layer perceptron (MLP) in which the output layer
is equipped with a non-linear activation function -- i.e. the RB solver --
in the sense that it does not contain any trainable parameters. Contrarily to standard activation functions, however,
the RB solver maps vectors -- which in our case contain either the physical
parameters of the PDE or the coefficients of the affine decomposition -- to (typically larger) vectors, i.e.
the values of the solution at the output locations. Alternatively, PDE-DNNs can be
regarded as asymmetric autoencoders in which the MLP and the RB solver play the roles
of the encoder and the decoder, respectively. In this view, the \textit{latent space} coincides with
either the space of physical parameters or of the space of coefficients of the affine decomposition,
depending on the structure of the underlying PDE.

This paper is structured as follows. In Section~\ref{sec:deepnns}, we recall basic concepts of NNs and
set the notation and the terminology that will be extensively adopted throughout
the rest of the paper. In Section~\ref{sec:parmpdes}, we provide a (non-exhaustive, but sufficiently
detailed for our purposes) introduction to the reduced order modeling of parametrized PDEs through the RB method.
In Section~\ref{sec:pde-aware}, we introduce the concept of PDE-DNN. In Section~\ref{sec:numexp}, we present numerical
experiments on a three-dimensional example of advection-diffusion allowing for an affine
decomposition (Section~\ref{subsec:affine}) and on two-dimentional non-affine examples, namely
a diffusion problem with a non-affine diffusion coefficient (Section~\ref{subsec:nonaffine}) and
the steady Navier-Stokes equations in a bifurcation with a resistive term modelling a stenosis (Section~\ref{sec:navier_stokes}).
Finally, in Section~\ref{sec:conclusions} some conclusions are drawn.

\section{Deep Neural Networks}
\label{sec:deepnns}
Artificial neural networks, often simply denoted neural networks (NNs) in the literature, are biologically inspired ML algorithms which are able to {\it learn} from observational data \cite{haykin2004neural}. Similarly to other ML techniques, NN models are described by means of a collection of trainable parameters. These parameters are inferred during a process whose goal is to find a locally (in the parameter space) optimal choice that leads to an ``accurate'' input-output mapping of the network on a dataset -- usually called training dataset -- for which the expected result is known.

\def\layersep{2.cm}
\begin{figure}
	\centering
	\begin{minipage}{0.65\textwidth}
	\def\layersep{2.cm}
	\begin{center}
		\resizebox{\textwidth}{!}{%
			\begin{tikzpicture}[shorten >=1pt,->,draw=black!50, node distance=\layersep]
			\tikzstyle{every pin edge}=[<-,shorten <=1pt]
			\tikzstyle{neuron}=[circle,fill=black!25,minimum size=17pt,inner sep=0pt]
			\tikzstyle{input neuron}=[neuron, fill=red!50];
			\tikzstyle{output neuron}=[neuron, fill=green!50];
			\tikzstyle{hidden neuron}=[neuron, fill=blue!50];
			\tikzstyle{annot} = [text width=4em, text centered]
			
			% Draw the input layer nodes
			\foreach \name / \y in {1,...,3}
			% This is the same as writing \foreach \name / \y in {1/1,2/2,3/3,4/4}
			\node[input neuron, pin=left:\Input{\y}] (I-\name) at (0,-\y cm) {};
			
			% Draw the hidden layer nodes
			\foreach \name / \y in {1,...,4}
			\path[yshift=0.5cm]
			node[hidden neuron] (H-\name) at (\layersep,-\y cm) {};
			
			% Draw the output layer node
			%\node[output neuron,pin={[pin edge={->}]right:Output}, right of=H-3] (O) {};
			\foreach \name / \y in {1,...,2}
			% This is the same as writing \foreach \name / \y in {1/1,2/2,3/3,4/4}
			%\path[yshift=0.5cm]
			\node[output neuron, pin={[pin edge={->}]right:\Output{\y}}, right of=H-3] (O-\name) at (\layersep, -0.5cm -\y cm) {};
			
			% Connect every node in the input layer with every node in the
			% hidden layer.
			\foreach \source in {1,...,3}
			\foreach \dest in {1,...,4}
			\path (I-\source) edge (H-\dest);
			
			% Connect every node in the hidden layer with the output layer
			%\foreach \source in {1,...,5}
			%\path (H-\source) edge (O);
			
			\foreach \source in {1,...,4}
			\foreach \dest in {1,...,2}
			\path (H-\source) edge (O-\dest);
			
			% Annotate the layers
			\node[annot,above of=H-1, node distance=1cm] (hl) {Hidden layer};
			\node[annot,left of=hl] {Input layer};
			\node[annot,right of=hl] {Output layer};
			\end{tikzpicture}
		}
	\end{center}
	
	\let\layersep\undefined

	\end{minipage}
	\caption{Example of DNN with one fully connected hidden layer}
	\label{fig:ex1hidden}
\end{figure}
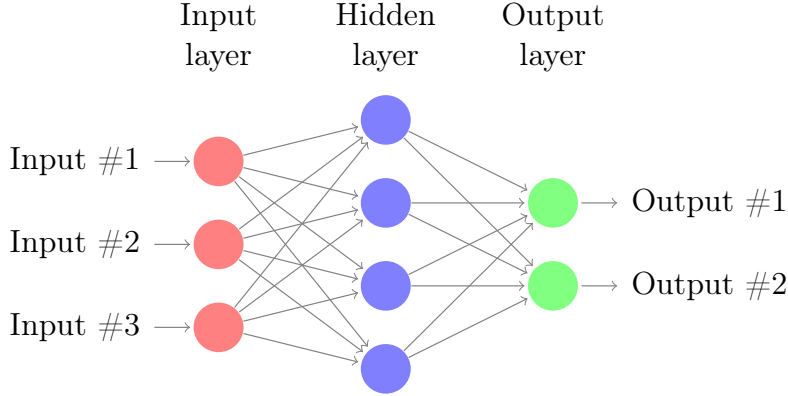

Deep Neural Networks (DNNs) are a class of NNs in which simple nonlinear modules are composed together in multiple layers. The first layer is the {\it input layer} and processes the raw data which feeds the network, while the last layer is the {\it output layer} which provides the output of the network; the ones located between the two are called {\it hidden layers}, see Figure~\ref{fig:ex1hidden}. A common choice for the architecture of the DNN is given by ordered fully connected layers of perceptrons, i.e. multi-layer perceptrons (MLPs). Perceptrons are simple computational units containing a set of nodes (the neurons) processing the data coming from the previous layer, and a non-linear {\it activation function}. MLPs are often preferred over other architectures for statistical regression, because of their ability to approximate non-trivial functions. It has been shown that MLPs with at least one hidden layer and differentiable activation functions are universal function approximators; in particular, MLPs with one hidden layer are able to approximate any continuous function, and MLPs with two hidden layers or more are able to approximate any function \cite{university1988continuous, cybenko1989approximation}.

Let us denote with $L$ the number of trainable layers in a MLP (i.e. to the hidden and output layers) and let us use the index 0 to chatacterize the input layer. The input layer performs a trivial mapping of the input of the neural network $\mathbf x^{(0)} \in \real^{\Nin}$ into itself.
The $l^\text{th}$ layer of a MLP, with $l > 0$, is a mapping of the form
\begin{equation}\label{eq:perceptron}
\yyl = \sigmal( \wl \xxl + \bl),
\end{equation}
where $\xxl \in \realNinl$ is the input, $\yyl \in \realNoutl$ is the output, $\wl \in \realNoutlNinl$ is a matrix of parameters (weights), $\bl \in \realNoutl$ is a vector of parameters (biases), and $\sigmal: \realNoutl \rightarrow \realNoutl$ is a non-linear activation function; note that $\xxl = \yylm$. In this work, we restrict ourselves to two of the most commonly used activation functions: $\relu(x) = \text{max}(0,x)$ and $\sigm(x) = (1 + \exp(-x))^{-1}$; applied to multidimensional inputs, these must be intended as component-wise functions. Equation~\ref{eq:perceptron} can be trivially extended to the case in which $\Ns$ samples are processed at the same time; in such scenario, the input and output of the layer are matrices -- that is $\XXl \in \real^{\Ns \times \Ninl}$ and $\YYl \in \real^{\Ns \times \Noutl}$ -- where each sample is stored row-wise.

We indicate $\thetal = (\wl,\bl)$ the weights-biases pair of the $l^\text{th}$ layer, and $\Theta = \{ \theta^{(1)},\,\ldots,\, \theta^{(L)} \}$ the set of all the parameters of a MLP.  Moreover, we denote $\Nout = \Nout^{(L)}$ the size of the output of the neural network. The goal of the training process for a MLP is to find suitable values of the parameters $\wl$ and $\bl$ for $l = 1,\ldots,L$  such that the function modeled by the network, $f: \realNin \rightarrow \realNout$ defined as $f(\mathbf x^{(0)}) := \mathbf y^{(L)}$, well approximates the actual input-output relationship.
During the training process, the discrepancy between the real output $\mathbf y$ and the one of the network $\mathbf y^{(L)}$ is measured by a {\it loss function} $\loss_\Theta(\mathbf y, \mathbf y^{(L)})$, which clearly depends on the parameters and is minimized during the training process. A common choice of loss function, which we also consider in this work, is given by the mean squared error (MSE). In the case of multiple samples, this is defined as
\begin{equation}
\MSE_\Theta(\mathbf Y, \mathbf Y^{(L)}) = \sum_{i = 1}^{\Ns} \sum_{j = 1}^{\Nout} \left ( Y_{ij} - Y^{(L)}_{ij} \right )^2,
\label{eq:mseloss}
\end{equation}
where $\mathbf Y$ is the matrix storing the desired outputs for all the samples. The minimization is typically performed via (stochastic) gradient descent or variations such as the Adam algorithm \cite{kingma2014adam}, which all require the computation of the gradients of the loss function with respect to the parameters. In modern DL frameworks, these quantities are automatically computed by {\it backpropagation}, which is a specialization of the chain rule for derivatives to the functionals modeled by such architectures.

\section{Parametrized partial differential equations}
\label{sec:parmpdes}
In this section we introduce the general framework of parametrized PDEs and we describe the fundamental principles of the reduced basis method.

\subsection{Parametrized PDEs}
\label{subsec:parametrized}
Let us consider a parameter space  $ \paramspace \subset \real^p $, $p \geq 1 $, and denote by $\param \in \paramspace $ a parameter vector encoding physical and/or geometrical properties of the problem. Furthermore, let us introduce an open and bounded domain $\Omega\subset \real^d, \,d =2,3 $, and denote by $\partial \Omega $ its boundary.
We  consider the following equation
\begin{align}\label{eq:nonlinearModel}
\nonlinear[\uu; \param] = f(\param) \qquad & \text{in } \Omega
\end{align}
where $\nonlinear[\cdot; \param]$ is a $\param $-differential operator which models a physical phenomenon. For any $\param \in \paramspace $, we are interested in the solution $\uu = \umu $.
Equation \eqref{eq:nonlinearModel} is provided with proper boundary conditions on $\partial\Omega$.

A numerical approximation to the solution $\uu $ of \eqref{eq:nonlinearModel} can be obtained by means of a full order model (FOM). The latter is derived, e.\ g.,  starting from the variational formulation of \eqref{eq:nonlinearModel} and employing a (Petrov-)Galerkin projection onto a finite-dimensional space. Examples of such FOMs are the spectral element method, the finite volume method and the finite element (FE) method. In our numerical experiments we consider the FE approximation as our reference {\it true} solution.
Solving the FE problem is equivalent to solving a nonlinear algebraic system of (generally) large dimension
\begin{align}\label{eq:nonlinearSystem}
\nonlinearMatrixmu\vecuh = \vecfhmu,
\end{align}
where $ \nonlinearMatrixmu \in \realhh$ is the parametrized nonlinear matrix, $\vecfh(\param) \in \realh $ the right hand side and $\vecuh = \vecuh(\param) \in \realh $ is the FE vector containing the value of the solution at the mesh nodes. The dimension $\Nh $ is in general very large when facing real-world applications, and can easily reach millions.
When facing parameter dependent problem, i.e.\ when the solution of the PDE \eqref{eq:nonlinearModel} is sought for many instances of the physical parameter $\param $, the computational load entailed by repeatedly solving the FE problem \eqref{eq:nonlinearSystem} can be unbearable. This issue can be overcome by Reduced Order Modelling (ROM) techniques, which allow to exploit the parameter dependence to boost the solution of the parametrized-problem at hand.

\subsection{The reduced basis method for parametrized PDEs}
\label{subsec:rbmethod}
Among all the ROM approaches, the RB method has been one of the most studied and exploited in the field of parametrized PDEs, and represents a convenient framework for the modeling reduction of differential problems.
In the following, we briefly recall the basic ideas of the RB method, we refer to \cite{quarteroni2015reduced,hesthaven2015certified} for full reviews and additional details.

The RB method relies on the idea that the $\param$-dependent solution of the $ \Nh \times \Nh $ FE problem can be well approximated by a linear combination of $N$ basis functions obtained as solution of the same problem for (suitably chosen) parameter values, with $N \ll \Nh$.
The computation is usually based on an {\it offline/online} splitting. The offline phase computes a {\it reduced space} (the RB space),
which is algebraically represented by the matrix $\rbV \in \real^{N_h\times N}$,  $ \rbV = [\rbvecbasis_1 | \dots | \rbvecbasis_N]$.
The goal of the RB method is to find a vector of weights $\rbvecu \in \realN$ such that $\rbV\rbvecu $ is an accurate approximation to the FE solution, that is $\vecuh \approx \rbV\rbvecu$.
This is done by solving an empirical algebraic RB problem of (small) dimension $N$ for any new instance of the parameter $\param$
\begin{align}\label{eq:rbAlgebraic}
\rbNonlinearMatrixmu\rbvecu = \rbvecfmu,
\end{align}
where  $\rbNonlinearMatrix \in \realNN $ is the RB matrix  and $\rbvecfmu \in \realN $ the RB right hand side, respectively obtained as  Galerkin projection of the orginal FE arrays, where the nonlinear operator is computed at the approximation $\rbV\rbvecu$
\begin{align}\label{eq:rbArrays}
\rbNonlinearMatrixmu = \rbV^T\nonlinearMatrix[\rbV\rbvecu; \param]\rbV,
\qquad
\rbvecfmu = \rbV^T\vecfhmu.
\end{align}
We notice that Petrov-Galerkin projections are also viable ways to obtain a robust RB approximation, especially for noncoercive problems, see e.g. \cite{abdulle2015petrov,DALSANTO2019186}.

\subsubsection{Constructing the RB space}

The RB matrix $\rbV$ can be obtained by employing either a greedy or proper orthogonal decomposition (POD) technique. We briefly remind the latter as we will employ it in the numerical experiments, for further reading we refer to \cite{volkwein2013proper}.
Let us consider a set of $\ns$ parameter values $\{\param_i\}_{i=1}^{\ns}$ and FE vectors $\setsnapi \subset \realh $ (called {\it snapshots}) collected as columns of a matrix $\snap = [\snaps_1| \dots |\snaps_{\ns}] $, $\snap \in \real^{\Nh \times \ns }$.
For any prescribed dimension $N$, the POD allows to find an orthonormal basis $\{\rbvecbasis_i\}_{i=1}^{\ns} $ and the corresponding $N$-dimensional subspace, spanned by the columns of the matrix $\rbV = [\rbvecbasis_1|\dots|\rbvecbasis_N] $ which best approximates $\setsnapi $ up to a tolerance $\epod$ with respect to the Euclidean norm. POD takes advantage of the singular value decomposition (SVD) of $\snap$
\begin{align*}
\snap = \svdU\svdSigma\svdZ^T,
\end{align*}
with $\svdU \in \realhh $ and $\svdZ \in \realnsns $ orthogonal matrices and $\svdSigma = diag(\sigma_1, \dots \sigma_{\ns}) $, $\svdSigma \in \realhns $, containing the singular values $\sigma_1 \geq \sigma_2 \geq \dots \geq \sigma_{\ns} \geq 0 $. Then,  $\rbV$ is provided by the first $N$ columns of $\svdU$, which form by construction an orthonormal basis for the best $N$-dimensional approximation subspace. In practical applications, $N$ is usually not set a priori by the user but for a fixed tolerance $\epod$, $N$ is found as the minimum $j$ such that
\begin{equation}
1 - \dfrac{\sum_{i = 1}^j \sigma_i^2}{\sum_{i = 1}^{\ns} \sigma_i^2} \leq \epod^2.
\end{equation}

\subsubsection{Affine decompositions of RB arrays}

It is clear that the solution of the RB problem \eqref{eq:rbAlgebraic} yields a number of operations which depends on the number of RB degrees of freedom $N$, hence independent of the FE dimension $\Nh$. However, the assembly \eqref{eq:rbArrays} of problem \eqref{eq:rbAlgebraic} requires the projection of the FE matrix and right hand side on the RB space, thus entailing operations with a complexity dependent on $\Nh $ and preventing an efficient {\it offline/online} splitting.
To overcome this bottleneck, the RB method relies on the affine dependence of the RB arrays, that is $\nonlinearMatrixmu $ and $\vecfmu $ can be written as sum of $\Qn $ and $\Qf $ terms independent of $\rbvecu $ and $ \param $, as follows
\begin{align}\label{eq:affinedependence}
\nonlinearMatrixuNmu =
\sum\limits_{q=1}^{\Qn} \thetanmu\nonlinearMatrixq,
\qquad
\vecfmu =
\sum\limits_{q=1}^{\Qf} \thetafmu\vecfq,
\end{align}
with $\nonlinearMatrixq \in \realhh, \, q = 1, \dots, \Qn $ and $\vecfq\in \realh, \, q = 1, \dots, \Qf $.
From \eqref{eq:affinedependence} it easily follows that
\begin{align}\label{eq:rbaffinedependence}
\rbNonlinearMatrixmu =
\sum\limits_{q=1}^{\Qn} \thetanmu\rbNonlinearMatrixq,
\qquad
\rbvecfmu =
\sum\limits_{q=1}^{\Qf} \thetafmu\rbvecfq,
\end{align}
with $\rbNonlinearMatrixq = \rbV^T \nonlinearMatrixq \rbV \in \realNN, \, q = 1, \dots, \Qn  $, and $\rbvecfq = \rbV^T\vecfq, \, q = 1, \dots, \Qf$.
Being the RB affine arrays  in \eqref{eq:affinedependence} $\rbvecu$- and $\param$-independent, they can be preassembled in the offline phase. Then the assembly of $\rbNonlinearMatrixmu $ and $\rbvecfmu $ in the online phase entails only the sums in \eqref{eq:rbaffinedependence}, whose operations are $\Nh$-independent.

When we are facing a PDE problem which does not feature by construction an affine decomposition as in \eqref{eq:affinedependence}, we cannot directly split the computation in an {\it offline} and {\it online} phase. However, an approximated affine decomposition can be computed by using the Empirical Interpolation Method (EIM) \cite{barrault2004empirical} or its discrete variants Discrete EIM (DEIM) and Matrix DEIM (MDEIM), where the latter is specific for matrices \cite{Chaturantabut2010deim,negri2015efficient}.
During the offline phase, these algorithms are provided with a set of vector (DEIM) or matrix (MDEIM) snapshots and a tolerance which encodes the accuracy of the resulting affine approximation and return a basis of affine components which is computed through an internal POD. Then, during the online phase an interpolation problem is solved to compute the coefficients $\thetanmu $ and $ \thetafmu$ for each new instance of the PDE problem.

Depending on the problem at hand and the chosen accuracy, the number of affine components can largely vary, from tens to hundreds or even thousands, which may largely affect the speedup provided by the RB approximation compared with the FE one.
This issue is prominent when facing nonlinear unsteady problems where the time-dynamics plays a relevant role for the creation of the RB approximation, see e.g. \cite{dalsanto2018hyper} in the case of the Navier-Stokes equations, where the number of affine components is large even for a relatively modest Reynolds number (about 400).

(M)(D)EIM techniques provide a satisfactory tool to deal with nonlinear and nonaffine PDE problems, however they present two significant bottlenecks:
\begin{itemize}
	\item the underlying FE mesh must always be available both in the offline and online phase; in particular, in the latter the mesh is essential for assembling the right hand side of the interpolation problem for computing $\thetanmu$ and $\thetafmu$. When the RB method is used to speedup the solution of large FE problems, the storage of the mesh can be significantly demanding and can prevent from using the RB approximation outside of a HPC environment.
	\item Given the (M)(D)EIM basis, the interpolation problem solved to produce  the affine decomposition does not necessarily provides the best combination of basis functions, and in practice, for complex problems, this happens only if a significantly large number of basis is employed.
\end{itemize}

\section{PDE-aware deep neural networks}
\label{sec:pde-aware}
In this section we introduce an original framework for approximating parametrized PDEs by coupling a DNN and a RB solver.

We focus on the following class of problems. Let us consider a physical system, mathematically modeled by a parametrized PDE as that in Equation~\eqref{eq:nonlinearModel}, for which we know the value of the solution $\uu$ at certain points $\Psetin = \{\pp_i\}_{i = 1}^{\Nin}$; in realistic scenarios, these measurements could be obtained from multiple sensors positioned in $\Omega$ or on its boundary $\partial \Omega $. In this work, we restrict ourselves to steady physical phenomena, since unsteady cases need further developments. The PDE is parametrized by  a vector of parameters $\param \in \paramspace \subset \real^p $, $p \geq 1$ (following the notation introduced in Section~\ref{subsec:parametrized}). We assume that the values of the parameters are not known, or that their direct measurements is either complex, expensive or invasive. We want to define a neural network which, starting from measurements, is able not only to predict values on different parts of the domain, but also to provide the related parameter vector $ \param$ or the coefficients of the affine decomposition of the reduced system in Equation~\eqref{eq:rbaffinedependence}. The train samples are generated by varying the parameters $\param$, whose value is not included in the dataset.
We are interested in predicting the values of the solution at a set of points $\Psetout = \{\ppo_i\}_{i = 1}^{\Nout}$, given $\uu(\pp_i;\param)$, $\pp_i \in \Psetin$ for $i = 1,\ldots,\Nin$, and at the same time to infer the value of the underlying parameter vector $\param$. In our numerical applications, we will explore the cases in which $\Psetin \cap \Psetout = \emptyset$ and $\Psetin \equiv \Psetout$.

We attempt the solution of such problems with PDE-aware deep neural networks (PDE-DNNs). The key idea of PDE-DNNs lies on the use of a PDE solver as building block of a DNN. This is motivated by the fact that standard DNNs employed for solving PDE problems do not exploit the underlying physics, and, to the best of our knowledge, the PDE plays a relevant role only in the definition of the loss function, as in \cite{raissiI2017physics,raissiII2017physics}. In particular, our proposition is to consider a MLP with $L$ trainable layers in which the output layer encodes the discrete mathematical model.

The first $L-1$ trainable layers of the network are mappings of the form of Equation~\eqref{eq:perceptron} in which $\sigmal = \relu$ for $l = 1,\ldots,L-1$. The mapping of the output layer takes the form
\begin{equation}\label{eq:lastlayer}
\yyL = \sigmaRB( \wL \xxL + \bL),
\end{equation}
where $\sigma_\text{RB}: \real^s \rightarrow \real^{\Nout}$ represents the function of a RB solver acting from a $s$-dimensional representation of the solution space to the space of the solution values at the output locations $\Psetout$. The RB solver $\sigma_\text{RB}$ constructs such mapping as follows. Let us assume without loss of generality that $s = \Qn + \Qf$ and that $\latent = [\thetann{1},\ldots,\thetann{\Qn},\thetaff{1},\ldots,\thetaff{\Qf}]\in \real^s$ is the vector of the coefficients of the affine decomposition in Equation~\eqref{eq:rbaffinedependence}. As a matter of fact, in the case the PDE allow for an affine decomposition, $\latent$ could contain the physical parameters themselves. Given $\latent$, the reduced solution is simply obtained by solving the linear system efficiently assembled as in Equation~\eqref{eq:rbaffinedependence}; we remark that the matrices $\rbNonlinearMatrixq$ and the vectors $\rbvecfq$ appearing in Equation~\eqref{eq:rbaffinedependence} must be computed a priori in the offline phase. Then, the desired output is found by projecting the reduced solution onto the FE space, i.e. $\vecuh = \rbV \rbvecu$, and by interpolating the FE solution at the points of interest in $\Psetout$. In our tests, for simplicity, the input and the outputs are taken as values of the FE solution at the nodes of the mesh; since we employ standard Lagrangian basis functions, no actual interpolation of the solution is needed, and, eventually, we can write the PDE-DNN activation function as follows
\newcommand{\Rout}{\mathbf{R}_{\operatorname{out}} }
\begin{align}
\sigmaRB( \latent ) = \Rout^T \rbV \rbvecu(\latent) = \Rout^T\rbV \rbNonlinearMatrix^{-1}(\latent)\rbvecf(\latent),
\end{align}
where $\Rout \in \real^{ \Nh \times \Nout}$ restricts the FE representation of the RB approximation to the output locations.
%\begin{align*}
%	&\Ndnn = \sum\limits_{q=1}^{\Qn} \thetann{q} \rbNonlinearMatrixq \approx \rbNonlinearMatrixmu, \qquad \fdnn = \sum\limits_{q=1}^{\Qf} \thetaff{q} \rbvecfq \approx \rbvecfmu,
%\end{align*}

We remark that, as we observed in our numerical simulations, taking $\xi \in [0,1]^s$ considerably improves the convergence rate of the optimization through gradient descent. For this reason, in our numerical simulations we actually consider $\widetilde \sigma_{\text{RB}} = \sigmaRB \circ \,\sigm$ as activation function in the output layer. The choices of $\relu$ and $\sigm$ as activation functions for the hidden layers and as preprocessing of the RB solver in the output layer respectively has been determined adequate after empirical observation. We aknowledge, however, that other architectures could lead to different and possibly better results than the ones reported in this paper.

Figure~\ref{fig:ex1hidden_pde} shows the schematic representation of a PDE-DNN with two trainable layers ($L = 2$) and with $s = 2$.

\def\nInput{3}
\def\nHidden{4}
\def\nOutput{2}
\def\nPdeOutput{3}

\begin{figure}
	\centering
	\begin{minipage}{0.65\textwidth}
	\def\nInput{3}
	\def\nHidden{4}
	\def\nOutput{2}
	\def\nPdeOutput{3}
	
	\def\layersep{2.cm}

	\begin{center}
		\resizebox{\textwidth}{!}{%
			\begin{tikzpicture}[shorten >=1pt,->,draw=black!50, node distance=\layersep]
			\tikzstyle{every pin edge}=[<-,shorten <=1pt]
			\tikzstyle{neuron}=[circle,fill=black!25,minimum size=17pt,inner sep=0pt]
			\tikzstyle{input neuron}=[neuron, fill=red!50];
			\tikzstyle{output neuron}=[neuron, fill=green!50];
			\tikzstyle{hidden neuron}=[neuron, fill=blue!50];
			\tikzstyle{annot} = [text width=4em, text centered]
			\tikzstyle{pde neuron}=[neuron, fill=cyan!50];
			\tikzstyle{write neuron}=[neuron, fill=white];
			\tikzstyle{write neuron draw}=[neuron, fill=white, draw];
			\draw[draw, fill=yellow!10] (3.3,-3.2) rectangle (7.5,-0.75)       ;
			
			% Draw the input layer nodes
			\foreach \name / \y / \coord in {1/1/103,2/2/144,3/3/456}
			% This is the same as writing \foreach \name / \y in {1/1,2/2,3/3,4/4}
			\node[input neuron, pin=left:$\vecuhimu{\coord} $ ] (I-\name) at (0,-\y cm) {};
			%	\node[input neuron ] (I-\name) at (0,-\y cm) {$\vecuhimu{\coord} $};
			
			% Draw the hidden layer nodes
			\foreach \name / \y in {1,...,\nHidden}
			\path[yshift=0.5cm]
			node[hidden neuron] (H-\name) at (\layersep,-\y cm) {};
			
			% Draw the output layer node
			%\node[output neuron,pin={[pin edge={->}]right:Output}, right of=H-3] (O) {};
			\foreach \name / \y in {1,...,\nOutput}
			% This is the same as writing \foreach \name / \y in {1/1,2/2,3/3,4/4}
			\path[yshift=0.5cm]
			node[output neuron, right of=H-3] (O-\name) at (\layersep, -1cm -\y cm) {};
			
			\foreach \name / \y in {1,...,\nOutput}
			\path[yshift=0.5cm]
			node[write neuron draw, right of=O-3] (W-\name) at  (3.25, -1cm -\y cm) {$\Dmu{\y}$};
			
			\node[pde neuron, right of=W-3] (pdesolver) at (4.75cm,-2cm) {};
			
			\foreach \name / \y / \coord in {1/1/143,2/2/167,3/3/236}
			\node[write neuron, right of=O] (PDE-\name) at (6.25cm,-\y) { $\vecuhimu{\coord}$ };

			% Connect every node in the input layer with every node in the
			% hidden layer.
			\foreach \source in {1,...,\nInput}
			\foreach \dest in {1,...,\nHidden}
			\path (I-\source) edge (H-\dest);
			
			% Connect every node in the hidden layer with the output layer
			%\foreach \source in {1,...,5}
			%\path (H-\source) edge (O);
			
			\foreach \source in {1,...,\nHidden}
			\foreach \dest in {1,...,\nOutput}
			\path (H-\source) edge (O-\dest);
			
			\foreach \source in {1,...,\nOutput}
			\path (O-\source) edge (W-\source);
			
			\foreach \source in {1,...,\nOutput}
			\path (W-\source) edge (pdesolver);
			
			\foreach \source in {1,...,\nPdeOutput}
			\path (pdesolver) edge (PDE-\source);
			
			% Annotate the layers
			\node[annot,above of=H-1, node distance=1cm] (hl) {Hidden layer};
			\node[annot,left of=hl] {Input layer};
			%	\node[annot,right of=hl] (ol) {Hidden layer};
			%	\node[annot, above of=pdesolver, node distance=2.5cm ] {PDE solver};
			\node[annot,right of=hl, node distance=3.34cm]   {Output layer};
			\node[annot] at (6.76,-2.6) {$\sigmaRB$};
			\end{tikzpicture}
		}
	\end{center}
	
	\let\nInput\undefined
	\let\nHidden\undefined
	\let\nOutput\undefined
	\let\nPdeOutput\undefined
	
	\let\layersep\undefined

	\end{minipage}
	\caption{Example of PDE-DNN with one hidden layer in the MLP.}
	\label{fig:ex1hidden_pde}
\end{figure}
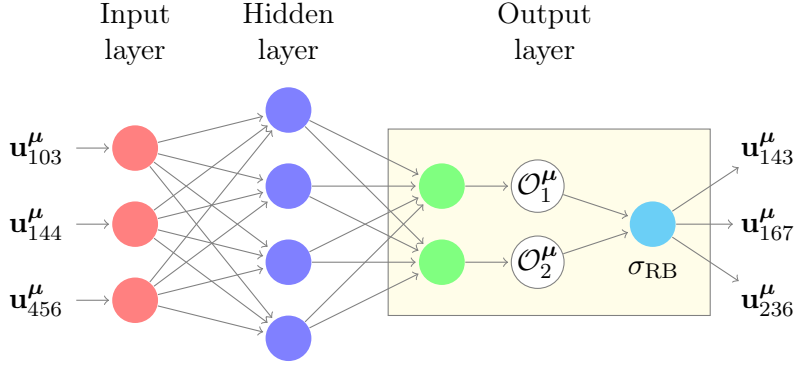

\begin{remark}
	Activation functions play a relevant role when training a network, since both the forward and backward propagation stages depend on them. 
	In particular, in the latter the derivatives of the activation function with respect to the latent space $\latent = [\thetann{1},\ldots,\thetann{\Qn},\thetaff{1},\ldots,\thetaff{\Qf}] $, that is $\frac{\partial \sigmaRB}{\partial \thetann{q}}$ and $\frac{\partial \sigmaRB}{\partial \thetaff{q}}$, enter into play. 
	Following the definition of the affine decompositions \eqref{eq:rbaffinedependence}, we can compute them as follows
	\begin{multline}\label{eq:backProp}
	\displaystyle
	\frac{\partial \sigmaRB}{\partial \thetann{q}}
	=
	\Rout^T \rbV \frac{\partial \Big(\rbNonlinearMatrix^{-1}(\latent)\rbvecf(\latent)\Big)}{\partial \thetann{q}}
	=
	\Rout^T \rbV \frac{\partial (\rbNonlinearMatrix^{-1}(\latent))}{\partial \thetann{q}}\rbvecf(\latent)
	\\
	=
	- \Rout^T \rbV \rbNonlinearMatrix^{-1}(\latent) \frac{\partial (\rbNonlinearMatrix(\latent))}{\partial \thetann{q}} \rbNonlinearMatrix^{-1}(\latent)\rbvecf
	=
	- \Rout^T \rbV \rbNonlinearMatrix^{-1}(\latent) \rbNonlinearMatrixq \rbNonlinearMatrix^{-1}(\latent)\rbvecf(\latent),
	\end{multline}
	for $q = 1, \dots, \Qn$, and, for $q = 1, \dots, \Qf$,
	\begin{equation*}
	\displaystyle
	\frac{\partial \sigmaRB}{\partial \thetaff{q}}
	=
	\Rout^T \rbV \frac{\partial \Big(\rbNonlinearMatrix^{-1}(\latent)\rbvecf(\latent)\Big)}{\partial \thetaff{q}}
	=
	\Rout^T \rbV \rbNonlinearMatrix^{-1}(\latent) \frac{\partial \Big(\rbvecf(\latent)\Big)}{\partial \thetaff{q}}
	=
	\Rout^T \rbV \rbNonlinearMatrix^{-1}(\latent) \rbvecfq.
	\end{equation*}
\end{remark}

Let us assume that we are able to build, either from direct measurements or -- as we consider in this paper -- from numerical simulations, a training dataset composed of $\Ns$ input vectors of the form $\xxi = [\uu(\pp_1;\param_i),\ldots,\uu(\pp_{\Nin};\param_i)]$ and $\Ns$ output vectors of the form $\yyi = [\uu(\ppo_1;\param_i),\ldots,\uu(\ppo_{\Nout};\param_i)]$, possibly arranged in matrix form in $\XX \in \real^{\Ns \times \Nin}$ and $\YY \in \real^{\Ns \times \Nout}$ respectively; the parameters $\param_i \in \paramspace$, for $i = 1,\ldots,\Ns$, are sampled in $\paramspace$ to well represent the whole set. In this regard, we note that in the case the training dataset be obtained from numerical simulation, the computational burden of the RB offline phase is partially mitigated by the fact that the numerical solutions corresponding to each $\param_i$ for $i = 1,\ldots,\Ns$ can be used as snapshots to build the RB basis by POD. For the training of the network, we consider the MSE loss function defined as in Equation~\eqref{eq:mseloss}.
%We note that, if the RB solver is designed so as to accept an estimate of $\param$ from the MLP, this choice of loss function (solely depending on the value of the function at a handful of points) makes our network able to estimate the value of the physical parameters without the need to provide such parameters in the training dataset.
For each sample, our loss function does not consider at all the values of the corresponding parameters $\param_i $. However, $\param_i$ or the coefficients of the affine decomposition are internally estimated by the MLP. Such additional information is extracted in our code by combining in the output tensor the solution by the RB solver with the result of the linear operation in the output layer (i.e. the input required by $\sigmaRB$). We stress that the latter is not used in the computation of the loss function.

\section{Numerical experiments}
\label{sec:numexp}
\subsection{Affinely parametrized elliptic problems}
\label{subsec:affine}
\begin{figure}
	\centering
	\includegraphics[scale=0.35,trim={0 3cm 0 3cm},clip]{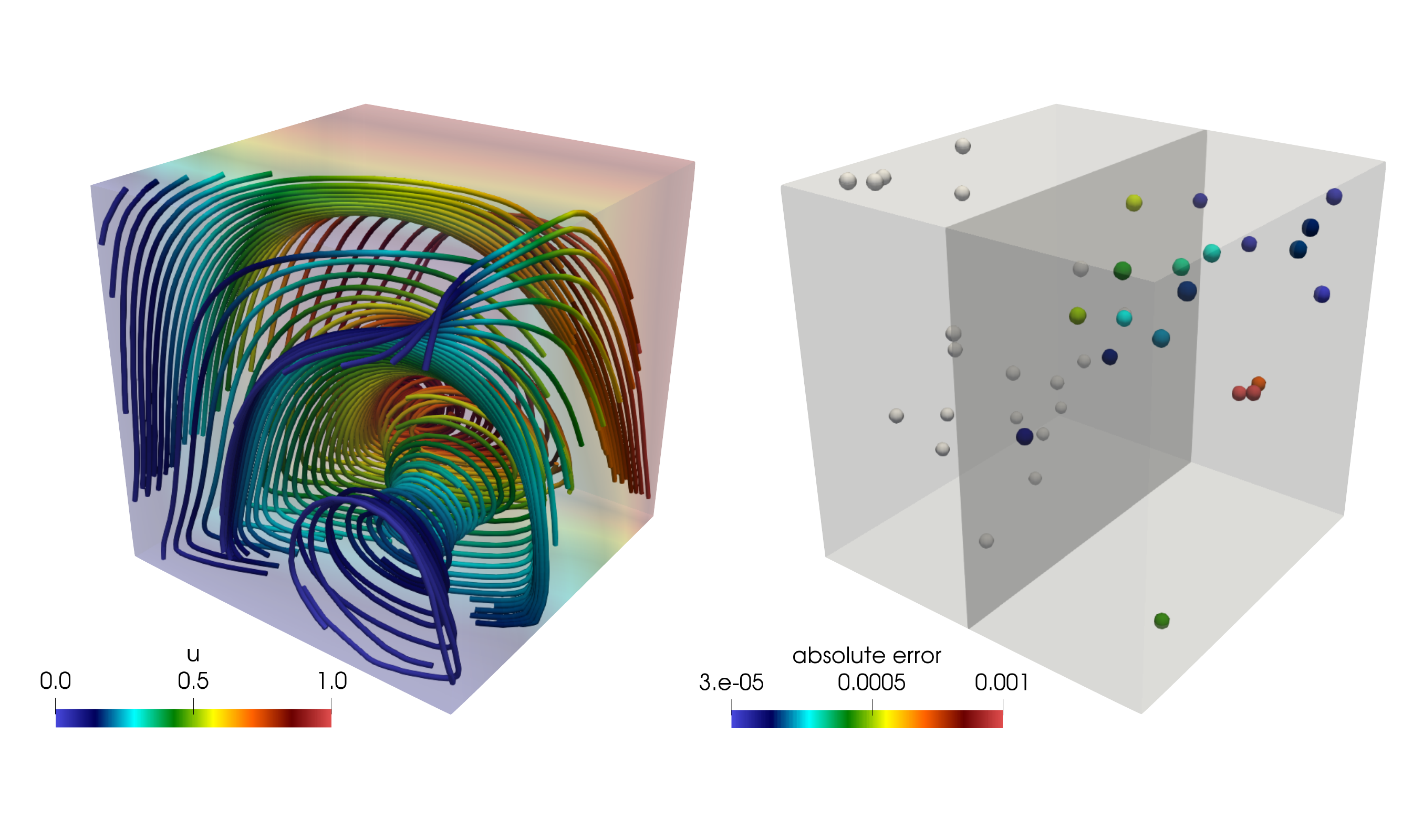}
	\caption{On the left, streamlines of the vector field $\vecb(\vecx;\alpha)$, colored with respect to the solution of Eq.~\ref{eq:testaffineadvection}. On the right, absolute error on the prediction of a PDE-DNN with 20 input points (white spheres) and 20 output points (colored spheres). Both figures refer to $\nu = 2.22$ and $\alpha = 0.48$.}
	\label{fig:advection_solution}
\end{figure}
Let us consider the domain $\Omega = (0,1)^3$ and the parametrized advection-diffusion problem
\begin{align}\label{eq:testaffineadvection}
\begin{cases}
-\nabla \cdot (\nu \nabla \uu) + \vecb(\vecx;\alpha) \cdot \nabla \uu = 0 \qquad & \text{ in } \Omega,\\
\mbox{+b.c.}
\end{cases}
\end{align}
with parameters $\param = (\nu, \alpha) \in \paramspace = (0.5, 10)\times(0, \pi/6) \subset \real^2$. The former defines the diffusivity of the problem and the latter is the angle defining the direction of the advection field $\vecb(\vecx;\alpha) = \sin(\alpha) \vecb_1(\vecx) + \cos(\alpha) \vec b_2(\vecx)$, which is obtained as linear combination of two divergence free $\alpha$-independent vector fields $\vecb_1(\vecx)$ and $\vecb_2(\vecx)$. These are computed as solutions of the Navier-Stokes equations with viscosity $1$ and equipped with non-homogeneous Dirichlet boundary conditions on one face of the cube -- specifically, $\vecb_1(\vecx) = [100,0,0]$ and $\vecb_2(\vecx) = [0,100,0]$ for every $\vecx \in \partial \Omega$ having $x_3 = 1$ -- and homogeneous Dirichlet boundary conditions everywhere else. Fig.~\ref{fig:advection_solution} (left) shows the streamlines of the vector field obtained by choosing $\alpha = 0.48$.
The boundary of the domain $\partial \Omega$ is partitioned in $\Gamma_D = \{\vecx \in \partial \Omega:\,x_1 = 0 \text{ or } x_1 = 1\}$ and $\Gamma_N$, such that $\partial \Omega = \Gamma_D \cup \Gamma_N$ with $\overset{\circ}{\Gamma}_D \cap \overset{\circ}{\Gamma}_N = \emptyset$. We impose non-homogeneous Dirichlet boundary conditions on $\Gamma_D$, in particular $\uu = 1$ for every $\vecx \in \Gamma_D$ with $x_1 = 1$ and $\uu = 0$ otherwise, and homogeneous Neumann conditions on $\Gamma_N$. Such problem could model, for instance, the temperature of a fluid in steady regime in which the values at the two Dirichlet boundaries are kept fixed and heat is freely exchanged at the Neumann boundaries. Fig.~\ref{fig:advection_solution} (left) shows the solution of Eq.~\eqref{eq:testaffineadvection} over the streamlines of the advecting vector field.

After the discretization of problem \ref{eq:testaffineadvection} by the FE method (which we consider for all of the high fidelity models presented in this paper), the system is rewritten in algebraic form as
\begin{equation}
\mathbf A \mathbf u = \mathbf F,
\end{equation}
where $\mathbf A = \nu \mathbf A_1 + \sin(\alpha) \mathbf A_2 + \cos(\alpha) \mathbf A_3 \in \real^{\Nh \times \Nh}$ is the affine decomposition of the finite element matrix; $\mathbf F \in \real^{\Nh}$ accounts for the boundary conditions and admits a similar affine decomposition.
\subsubsection{PDE-DNNs for affinely parametrized PDEs}
We aim at evaluating the performance of a PDE-DNN in predicting solutions of Eq.~\eqref{eq:testaffineadvection} for several instances of $ \param$ in a set of output points in terms of the following \textit{hyperparameters}: {\em{i)}} the number of input points $\Nin$, {\em{ii)}}  the number of output points $\Nout$, {\em{iii)}}  the tolerance for the POD of the RB solver $\epod$, and {\em{iv)}} the number of samples in the training dataset $\Ns$.

For this test case we randomly sample the input and output points in the partitions of $\Omega$ corresponding to $x_2 \leq 1/2$ and $x_2 > 1/2$ respectively, as depicted in Fig.~\ref{eq:testaffineadvection} (right); this choice is motivated by the fact that, if two sets $\Psetin$ and $\Psetout$ were obtained from sampling over the entire $\Omega$, then also simple interpolation between points in the input set could potentially provide a rough approximation of the solution in the output points. As mentioned in Section~\ref{sec:pde-aware}, the sampling is done over the nodes of the mesh which is used in the computations by taking care of not selecting nodes on the Dirichlet boundaries; we consider a computational mesh composed of tetrahedra with $\mathbb{P}^1 $ finite elements and a total of $\Nh = 12416$ degrees of freedom.

\begin{figure}
	\centering
	\includegraphics[scale=0.5,trim={0cm 0 0cm 0},clip]{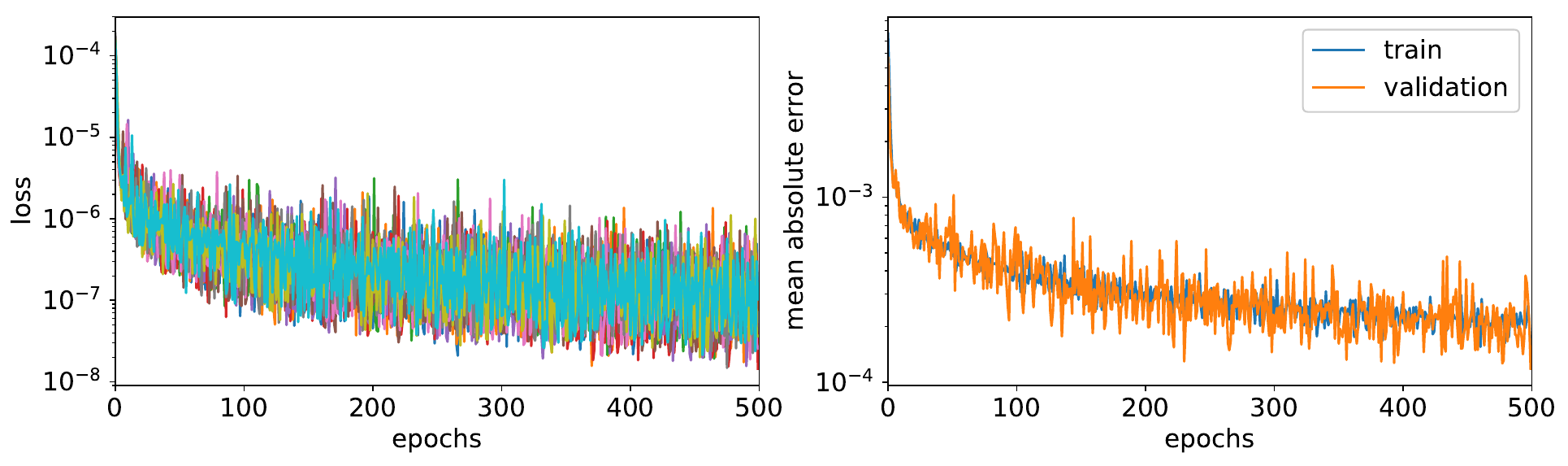}
	\caption{On the left, decaying of the MSE function during 10 trainings of PDE-DNN model with choice of hyperparameters $\Nin = 20$, $\Nout = 20$, $\epod = 10^{-5}$ and $\Ns = \text{10000}$. On average, the training has taken 1366 s. On the right, decaying of the mean absolute of error (averaged on the 10 runs) for train and validation datasets. The validation is obtained from the training dataset by selecting 20\% of the samples.}
	\label{fig:loss}
\end{figure}

In order to build a PDE-DNN, we first need to perform the offline phase of the RB method. This entails the collection of a predetermined number of snapshots $\ns$, which are obtained by solving the problem by FOM -- in our case, the FE method with linear Lagrangian basis functions -- for $\param_i \in \paramspace$, $i = 1,\ldots,\ns$. In this example, $n_s = 350$ and the parameters are sampled from a uniform distribution in $\paramspace$. Performing the POD of the snapshots matrix $\snap$ with tolerances $\epod = 10^{-4}$, $\epod = 10^{-5}$, $\epod = 10^{-6}$ and $\epod = 10^{-7}$ yields bases of sizes $N = 12$, $N = 19$, $N = 28$ and $N = 38$ respectively.

The $\ns = 350$ snapshots generated during the RB offline phase are also used as samples in the training dataset, which is performed through Adam algorithm with learning rate $10^{-3}$. We decide to enlarge such dataset by exploiting the ability of the RB method to generate solutions corresponding to new parameters in $\paramspace$ at a negligible cost. The largest training dataset we consider is composed of 20000 samples: 350 full order solutions and 19650 RB solutions obtained with POD tolerance $\epod = 10^{-7}$. As we wish to investigate the role of the number of samples on the training process and performance of the network, we also consider smaller training datasets by selecting the first $\Ns$ samples of the largest training dataset.
We remark that in our tests we actually select 20\% of the training samples for the validation of the model, that is, the training process is based on the remaining 80\% samples and the loss and other performance metrics are computed also on the validation dataset. In the deep learning community, such practice is adopted to prevent common issues such as \textit{over-} or \textit{underfitting}; we refer the reader to \cite{goodfellow2016deep} for more information on these topics. In the following, we will refer to training dataset as the subset of the total number of samples that is used during training. The test dataset is composed of 600 FE solutions at different random parameter values.

\begin{figure}
	\centering
	\includegraphics[scale=0.5,trim={0cm 0 0cm 0},clip]{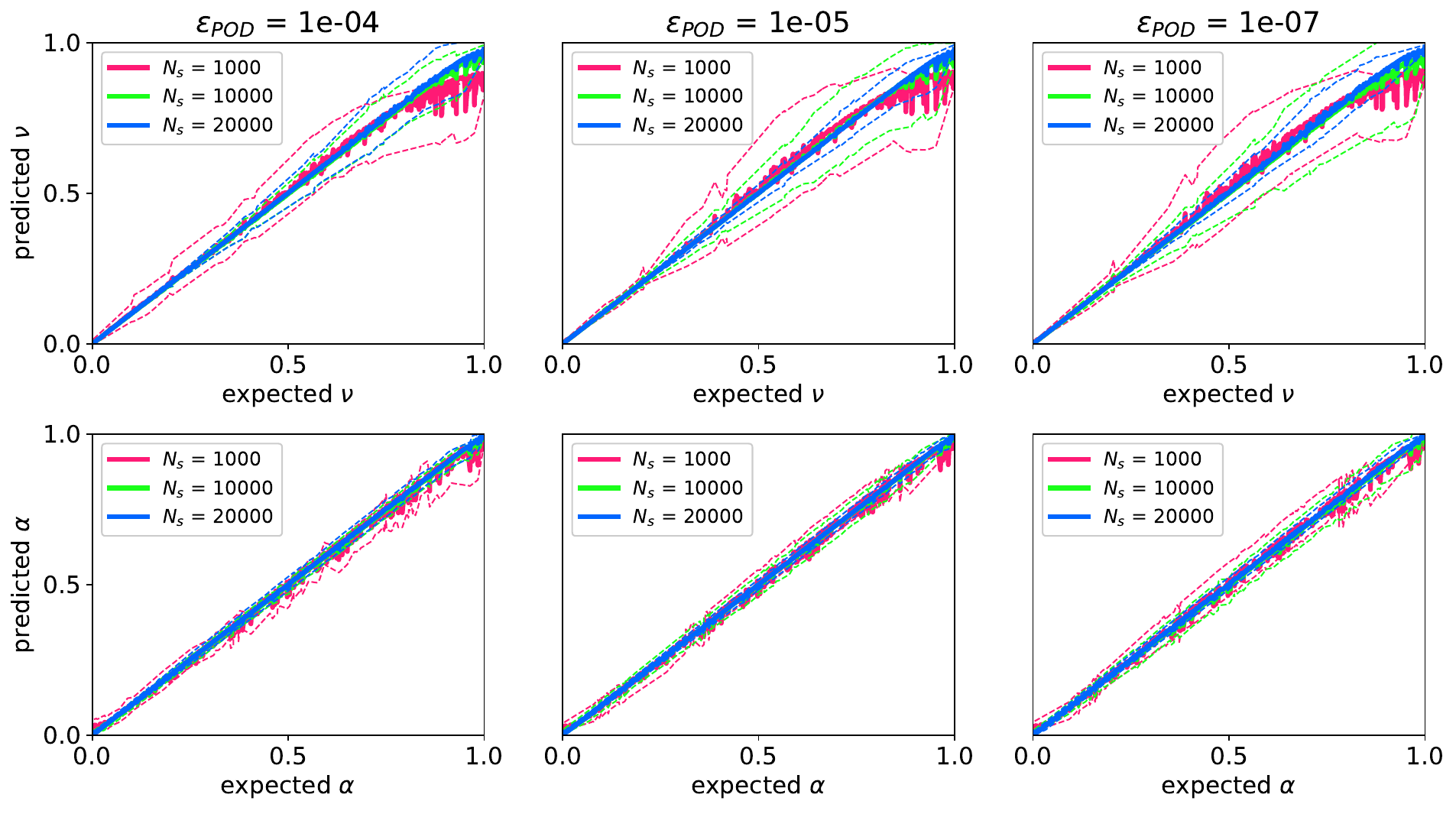}
	\caption{Estimation of the two (normalized over the respective ranges) physical parameters of the problem, $\nu$ (top row) and $\alpha$ (bottom row), for different POD tolerances of the RB problem and different number of samples. The results refer to the test dataset, which is composed of 600 data points, and are average over 30 trainings of the PDE-DNN network. The regions between dashed lines display the 95\% confidence intervals corresponding to the average lines with the matching colors.}
	\label{fig:parameter_estimation}
\end{figure}

The PDE-DNN consists of 4 hidden layers. In this test case, each hidden layer consists of 256 neurons. The linear part of the output layer is in charge of estimating a two dimensional vector, which is interpreted by the RB solver as an estimation of the physical parameters of the problem. The training of the PDE-DNN is performed through minimization of the loss in Eq.~\eqref{eq:mseloss} by gradient descent for 500 epochs; such number has been determined as adequate after empirical observations. Fig.~\ref{fig:loss} depicts the decaying of the MSE loss over 10 trainings of a PDE-DNN; we remark that the evolution of the loss function is not deterministic because the initialization of the trainable parameters of the MLP is random. In Fig.~\ref{fig:loss} (right), we focus on the mean absolute error -- averaged over the 10 runs -- as metric of the performance of the network. Since the mean absolute error reaches a plateau at around 500 epochs, we safely conclude that underfitting does not affect the training and this number of epochs is sufficient. Moreover, as the mean average error on the training and the validation dataset follow the same trend, overfitting does not occur during training. This is also true for the other choices of hyperparameters that we have considered. Our interpretation of this pheonomenon is that, as the train samples are obtained from numerical simulations, there are no
``hidden parameters'' that could characterize only the training dataset and penalize the performance over the validation. In other words, well representing the training data automatically leads to good approximation of the validation dataset.

Fig.~\ref{fig:parameter_estimation} displays the estimation of the normalized parameters $\nu$ and $\alpha$ performed by the PDE-NN for different values of the POD tolerance and different number of training samples. The findings are obtained by averaging the performance of 30  networks (trained over 500 epochs) on the test dataset. The effect of the two hyperparameters $\epod$ and $\Ns$ is the following. Given a sufficiently high number of input samples, for instance $\Ns = 20000$, the RB tolerance does not significantly affect the accuracy of the estimation for neither of the two physical parameters (the tolerance regions corresponding to $\Ns = 20000$, indeed, do not change significantly across the two rows of plots of Fig.~\ref{fig:parameter_estimation}). If, however, the training dataset is scarce, decreasing the tolerance of the POD leads to more uncertain results (i.e. larger confidence regions) without an appreciable increase of the precision on average. In other words, based on these results we conclude that the choice $\epod = 10^{-4}$ or $\epod = 10^{-5}$ is sufficient to get accurate results regardless of the size of the training dataset.

The estimation of the diffusion coefficient $\nu$ appears more challenging than the one of the angle $\alpha$: for values approaching $\nu = 10$  the networks are not able to distinguish accurately diffusivity coefficients close to each other, and in particular the predicted values are on average underestimating the correct ones (this is particularly true for $\Ns = 1000$). The estimation of the parameter $\alpha$ is quite accurate for all the considered configurations.

\begin{table}
	\centering
	\begin{tabular}{ c c c c c c c c c}
		\toprule
		& & \multicolumn{7}{c}{$(\Nin,\Nout)$} \\
		\cmidrule(l){3-9}& $\Ns$ & $(20,20)$ & $(40,40)$ & $(100,100)$ & $(20,40)$ & $(20,100)$ & $(40,20)$ & $(100,20)$ \\
		\midrule
		\multirow{3}{*}[-2pt]{PDE-NN} & 1000 & 2.18e-02 & 2.64e-02 & 2.00e-02 & 2.03e-02 & 2.31e-02 & 2.29e-02 & 7.05e-02 \\
		& 10000 & 3.62e-03 & 5.28e-03 & 5.00e-03 & 1.11e-02 & 9.99e-03 & 4.60e-03 & 3.87e-03 \\
		& 20000 & 2.87e-03 & 7.33e-03 & 5.69e-03 & 2.24e-03 & 4.23e-03 & 4.81e-03 & 3.50e-03 \\
		\midrule
		\multirow{3}{*}[-2pt]{MLP} & 1000 & 2.37e-02 & 1.91e-02 & 1.78e-02 & 2.24e-02 & 1.62e-02 & 1.75e-02 & 1.56e-02 \\
		& 10000 & 5.81e-03 & 7.49e-03 & 3.69e-03 & 6.84e-03 & 4.46e-03 & 4.70e-03 & 5.32e-03 \\
		& 20000 & 4.71e-03 & 3.96e-03 & 3.53e-03 & 4.17e-03 & 3.55e-03 & 4.07e-03 & 4.08e-03 \\
		\midrule
		\multirow{3}{*}[-2pt]{MLP$_\mu$} & 1000 & 2.39e-02 & 2.51e-02 & 3.04e-02 & 1.87e-02 & 2.01e-02 & 2.61e-02 & 2.02e-02 \\
		& 10000 & 7.99e-03 & 4.13e-03 & 5.80e-03 & 4.93e-03 & 6.50e-03 & 9.62e-03 & 6.12e-03 \\
		& 20000 & 2.97e-03 & 1.30e-02 & 2.82e-03 & 5.56e-03 & 8.53e-03 & 5.95e-03 & 3.17e-03 \\
		\bottomrule
	\end{tabular}
	\caption{Average error on $\nu$ over 600 test samples and over the outputs of 10 independently trained networks.}
	\label{tab:nu}
\end{table}
\begin{table}
	\begin{tabular}{ c c c c c c c c c}
		\toprule
		& & \multicolumn{7}{c}{$(\Nin,\Nout)$} \\
		\cmidrule(l){3-9}& $\Ns$ & $(20,20)$ & $(40,40)$ & $(100,100)$ & $(20,40)$ & $(20,100)$ & $(40,20)$ & $(100,20)$ \\
		\midrule
		\multirow{3}{*}[-2pt]{PDE-NN} & 1000 & 8.42e-03 & 1.58e-02 & 9.16e-03 & 9.52e-03 & 1.09e-02 & 1.30e-02 & 3.93e-02 \\
		& 10000 & 3.05e-03 & 2.36e-03 & 2.91e-03 & 4.71e-03 & 3.30e-03 & 1.69e-03 & 2.07e-03 \\
		& 20000 & 2.66e-03 & 4.28e-03 & 2.30e-03 & 1.02e-03 & 1.71e-03 & 4.54e-03 & 1.41e-03 \\
		\midrule
		\multirow{3}{*}[-2pt]{MLP} & 1000 & 8.44e-03 & 7.62e-03 & 1.46e-02 & 1.02e-02 & 8.10e-03 & 7.03e-03 & 8.61e-03 \\
		& 10000 & 7.10e-03 & 3.02e-03 & 4.49e-03 & 4.22e-03 & 1.81e-03 & 4.41e-03 & 2.43e-03 \\
		& 20000 & 2.31e-03 & 1.67e-03 & 1.34e-03 & 1.85e-03 & 2.46e-03 & 2.27e-03 & 1.85e-03 \\
		\midrule
		\multirow{3}{*}[-2pt]{MLP$_\mu$} & 1000 & 1.57e-02 & 1.07e-02 & 2.20e-02 & 8.79e-03 & 1.04e-02 & 1.21e-02 & 1.49e-02 \\
		& 10000 & 3.03e-03 & 1.95e-03 & 4.83e-03 & 2.08e-03 & 2.09e-03 & 2.67e-03 & 6.47e-03 \\
		& 20000 & 2.09e-03 & 3.64e-03 & 1.96e-03 & 2.11e-03 & 3.88e-03 & 4.26e-03 & 3.03e-03 \\
		\bottomrule
	\end{tabular}

	\caption{Average error on $\alpha$ over 600 test samples and over the outputs of 10 independently trained networks.}
	\label{tab:alpha}
\end{table}

%In this paper, we restrict ourselves to MLPs composed of 5 layers, in which $\Nin^{(0)} = \Nin$ and $\Nout^{(4)}$ is the number of values that we are interested in estimating.
Table \ref{tab:nu} and \ref{tab:alpha} show the errors on the prediction of $\nu$ and $\alpha$ respectively for different configurations of the PDE-DNN, averaged over 10 trainings of the network; all the results are obtained with $\epod = 10^{-5}$.

As a comparison, we consider three neural networks, denoted MLP$_{\mu}$, MLP$_{\operatorname{out}}$ and MLP. These are are designed and trained to compute the parameter $\param $ (MLP$_{\mu}$), or the  value of $\uu$ at the given output locations (MLP$_{\operatorname{out}}$) or both (MLP). The hidden layers of these networks are the same as those in the considered PDE-DNN, whereas the output layer is a perceptron equipped with $\sigm$ and with dimension  $\text{dim}(\boldsymbol \mu)$, $\Nout$ and $\Nout + \text{dim}(\boldsymbol \mu)$ for MLP$_{\mu}$, MLP$_{\operatorname{out}}$ and MLP, respectively.

As we noted in Fig.~\ref{fig:parameter_estimation}, the number of samples plays a relevant role in the approximation of the physical parameters, whereas the number of input and output locations does not affect significantly the results.
We remark that the PDE-DNN is able to obtain the approximately same accuracy of the other two neural networks, even though the parameters are only found as a byproduct (i.e. it is not necessary to train the model by providing the value of the parameters).

Table \ref{tab:output} shows the error on the output normalized with respect to the norm of the expected output (to take into account the different output sizes) in the same configurations considered in the previous tables.
Also in the case of MLP$_{\operatorname{out}}$, the three networks perform similarly for any choice of the hyperparameters.

\begin{table}
	\centering
	\begin{tabular}{ c c c c c c c c c}
		\toprule
		& & \multicolumn{7}{c}{$(\Nin,\Nout)$} \\
		\cmidrule(l){3-9}& $\Ns$ & $(20,20)$ & $(40,40)$ & $(100,100)$ & $(20,40)$ & $(20,100)$ & $(40,20)$ & $(100,20)$ \\
		\midrule
		\multirow{3}{*}[-2pt]{PDE-NN} & 1000 & 7.59e-04 & 8.72e-04 & 5.75e-04 & 5.74e-04 & 6.35e-04 & 1.57e-03 & 4.84e-03 \\
		& 10000 & 2.08e-04 & 2.70e-04 & 3.77e-04 & 3.18e-04 & 3.30e-04 & 1.80e-04 & 1.28e-04 \\
		& 20000 & 1.63e-04 & 2.21e-04 & 1.78e-04 & 7.88e-05 & 1.24e-04 & 2.66e-04 & 1.40e-04 \\
		\midrule
		\multirow{3}{*}[-2pt]{MLP} & 1000 & 1.54e-03 & 1.37e-03 & 2.56e-03 & 1.17e-03 & 9.85e-04 & 1.03e-03 & 1.41e-03 \\
		& 10000 & 8.93e-04 & 5.58e-04 & 5.34e-04 & 5.47e-04 & 3.48e-04 & 5.66e-04 & 4.34e-04 \\
		& 20000 & 2.98e-04 & 3.39e-04 & 2.89e-04 & 2.43e-04 & 2.28e-04 & 3.51e-04 & 3.93e-04 \\
		\midrule
		\multirow{3}{*}[-2pt]{MLP$_{\operatorname{out}}$} & 1000 & 1.05e-03 & 8.47e-04 & 1.45e-03 & 1.01e-03 & 1.34e-03 & 1.22e-03 & 1.12e-03 \\
		& 10000 & 5.99e-04 & 4.51e-04 & 1.81e-04 & 2.82e-04 & 1.56e-04 & 2.80e-04 & 1.99e-04 \\
		& 20000 & 1.63e-04 & 8.50e-05 & 1.78e-04 & 3.68e-04 & 1.33e-04 & 1.37e-04 & 1.11e-04 \\
		\bottomrule
	\end{tabular}
	\caption{Average error -- normalized with respect to the output norm -- on $\yyi$ over 600 test samples and over the outputs of 10 independently trained networks.}
	\label{tab:output}
\end{table}

\begin{figure}
	\centering
	\includegraphics[scale=0.57,trim={0cm 0 0cm 0},clip]{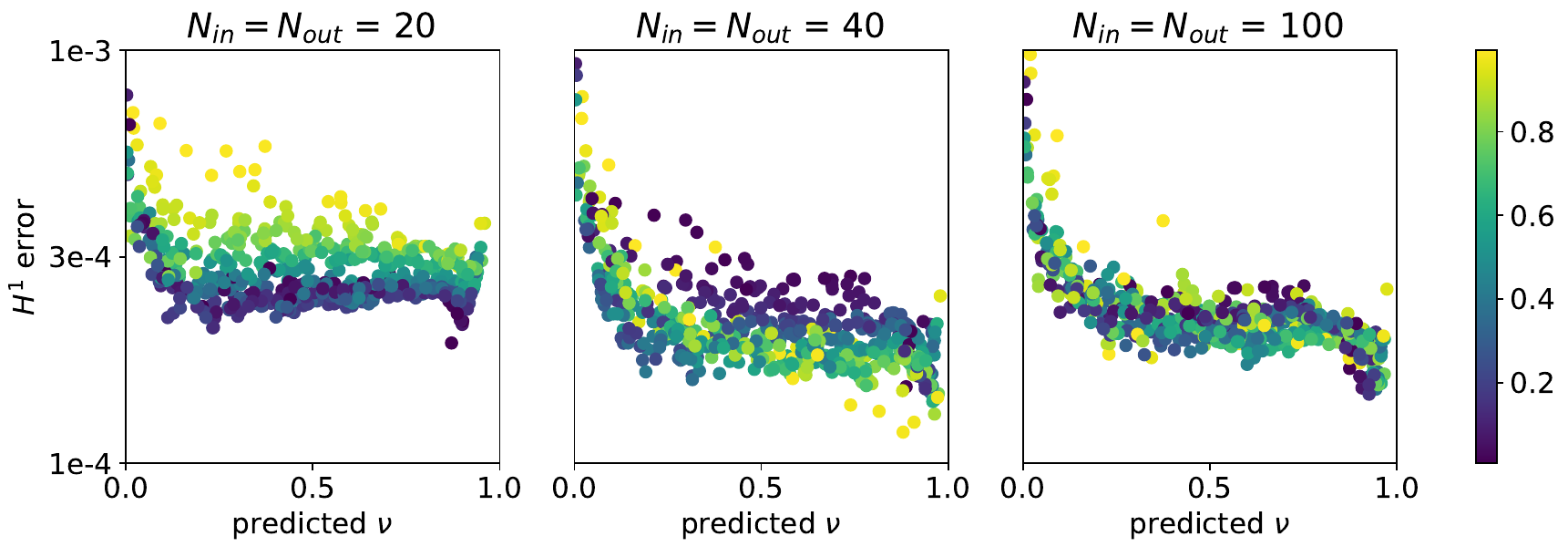}
	\caption{$H^1$ error against the 600 full order solutions in the test dataset when employing 20, 40 and 100 input locations in the dataset. Each point has $x$-component equal to the estimated value of the normalized physical parameter $\nu$ and it is colored with respect to the estimated value of the normalized physical parameter $\beta$; these and the errors have been obtained as averages over 10 trainings (500 epochs) of PDE-DNNs corresponding to $\epod = 10^{-5}$ and $Ns = 20000$.}
	\label{fig:h1error}
\end{figure}

\medskip
Finally, we consider the case in which the physical input and output location coincide. In other words, we train the network such that the RB solution coming from the output layer is as close as possible to the full order solution at the input locations($x\leq 1/2$). In this scenario, the PDE-NN is trained to model the identity function and the architecture is effectively an autoencoder in which the encoder is the MLP, the decoder is the RB solver, and the latent space is the space of physical parameters. In order to evaluate the performance of the network, we consider as metrics the $H^1$ error against the full order solution. We recall that this is possible because, given the estimation of the physical parameters by the MLP, the RB solver is able to approximate the solution at any point in the domain.

Fig.~\ref{fig:h1error} shows the $H^1$ error averaged over 10 trainings of PDE-DNNs when employing 20, 40 and 100 physical points. The average $H^1$ errors over the whole test dataset in the three cases is $3.0\times10^{-4}$, $2.3\times 10^{-4}$ and $2.4\times10^{-4}$. As a comparison, we report that the average error achieved with the RB method alone with the same POD accuracy $\epod = 10^{-5}$ is $1.1\times10^{-5}$. We observe that increasing the number of input locations has the effect of slightly decreasing the average error performed over the test dataset. Moreover, Fig.~\ref{fig:h1error} highlights that the prediction of the networks is considerably worse when the normalized $\nu$ approaches 0. This phenomenon and the fact that the loss function we consider is only concerned with the accuracy on the solution could explain the tendency of the networks to be considerably more precise in the estimation of $\nu$ when it is small (as displayed in Fig~\ref{fig:parameter_estimation}).

\subsection{Nonaffinely parametrized elliptic problems}
\label{subsec:nonaffine}
\begin{figure}
	\centering
	\subfloat{\includegraphics[width=0.45\textwidth]{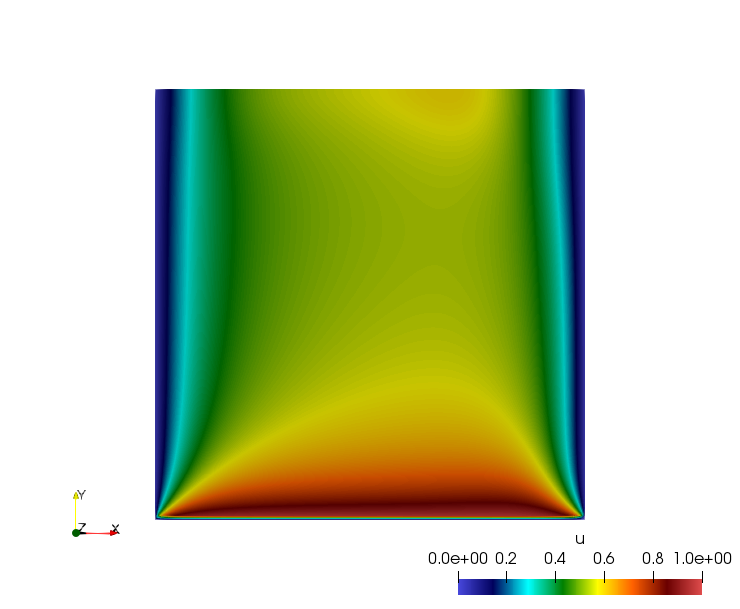}}
	\subfloat{\includegraphics[width=0.45\textwidth]{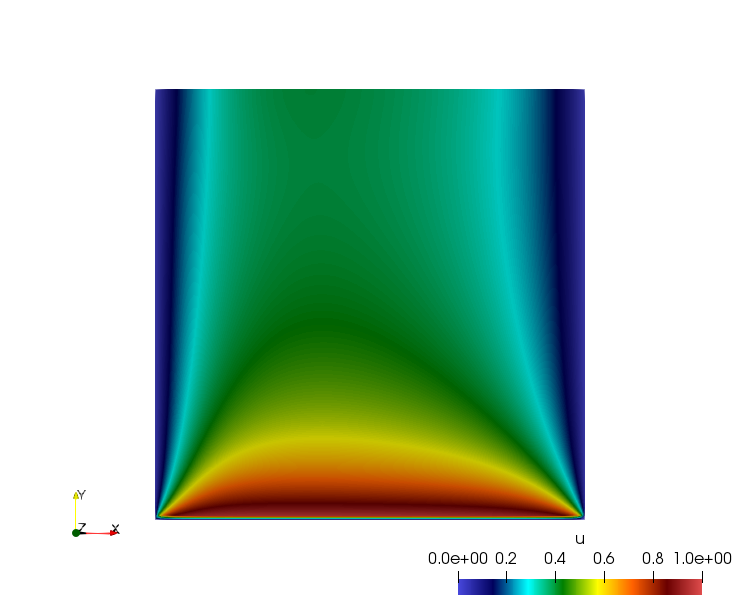}}
	\caption{Examples of solutions of problem \eqref{eq:testcaseII} for the parameters $\param_1 =(0.42, 0.42, 0.06) $ (left) and $\param_2 =(0.45, 0.59, 0.09) $ (right). }
	\label{fig:exmpleTestII}
\end{figure}

In the second test case, we consider a second-order diffusion problem where the diffusion coefficient is nonaffinely parametrized.
Let us consider the domain $\Omega = (0,1)^2 $ describing a square beam and a diffusion problem of the form
\begin{align}\label{eq:testcaseII}
\begin{cases}
-\nabla\cdot(\diffCoeffmu\nabla T) = 1 \qquad & \text{ in } \Omega, \\
%\diffCoeffmu\nabla T\cdot\vec{n} = 0       \qquad & \text{ on } \Gamma_N
%\\
\mbox{+b.c.}
%T = 0   \qquad & \text{ on } \Gamma_D^0,
%\\
%T = 1   \qquad & \text{ on } \Gamma_D^1,
\end{cases}
\end{align}
where $T$ is the temperature of the beam.
We define the following boundaries
\begin{align*}
\Gamma_N = \big\{\vecx \in \bar{\Omega}: \; y = 1\big\}, \qquad
\Gamma_D^1 = \big\{\vecx \in \bar{\Omega}: \; y = 0\big\}, \qquad
\Gamma_D^0 = \partial\Omega \backslash \Gamma_N \backslash \Gamma_D^1,
\end{align*}
such that $\partial\Omega = \bar{\Gamma}_0^D \cup \bar{\Gamma}_1^D \cup \bar{\Gamma}_N $, and we set homogeneous Neumann boundary conditions on $\Gamma_N$, homogeneous Dirichlet conditions on $\Gamma_D^0 $ and the solution to be equal to 1 on $\Gamma_D^1$.
We introduce the parameter vector
\begin{align*}
\param & = (x_0^1, x_0^2, \sigma) \in [0.4, 0.6]^2\times[0.05, 0.1] \subset \real^3
\end{align*}
and the coefficients
\begin{align*}
& \diffCoeffmu  = \sigma + \frac{1}{\sigma}\exp\left( -\frac{\|\vecx - \vecx_0(\param)\|^2}{\sigma} \right), \qquad \vecx_0(\param) = (x_0^1, x_0^2).
\end{align*}
Examples of solutions for two different physical parameters $\param$ are reported in Figure \ref{fig:exmpleTestII}.

We employ a lifting function $T_l = T_l(\vecx) $ to deal with the nonhomogeneous Dirichlet boundary conditions, and we use a structured mesh with $ \Nh$ = 10201 vertices and the FE method with first order polynomial piecewise basis functions to discretize the variational problem corresponding to \eqref{eq:testcaseII}. This leads to the following parametrized linear system of dimension  $\Nh $
\begin{align}\label{eq:nonaffineFem}
\feAhNamu \feTNa = \feFNamu,
\end{align}
where $\feAhNamu\in \realhh$ and $\feTNa=\feTNamu,\, \feFNamu\in \realh $.
A significant difference with respect to the previous test case is the nonaffine dependence of the FE matrix $\feAhNamu $ and right hand side $\feFNamu $ (the latter due to the lifting function) with respect to the parameters, i.e. assumptions \eqref{eq:affinedependence} do not hold. Hence, to build a RB model independent of dimension $\Nh $, we employ MDEIM and DEIM leading to $\QaNa $ sand $ \QfNa$ affine basis for the FE matrix and right hand side, respectively, that is
\begin{align}\label{eq:affine_approx_na}
\feAhNamu \approx \sum\limits_{q=1}^{\QaNa} \thetaAhNamu{q}\feAhNaq,
\qquad
\feFNamu \approx \sum\limits_{q=1}^{\QfNa} \thetaFNamu{q}\feFNaq.
\end{align}
We recall that given a new parameter $\param $, the coefficients $\{\thetaAhNamu{q}\}_{q=1}^{\QaNa} $ and $\{\thetaFNamu{q}\}_{q=1}^{\QfNa} $ are computed by solving an interpolation problem which enforces the true values of some preselected entries of the matrix (or right hand side) in its approximation. Clearly, the higher the number of affine components $\QaNa $ and $ \QfNa$, the more accurate the resulting affine approximation and the corresponding RB model. Notice that the underlying parametrization of $\diffCoeffmu$ results in a complex parameter dependence, especially related to the approximate affine approximation of $ \feAhNamu$. This is due to the Gaussian-like shape of the diffusion coefficient, which changes its centre and amplitude according to the value of the parameter $\param $ and requires MDEIM to compute a large number of basis functions to obtain an accurate affine approximation. This fact can hamper the online efficiency of the RB method and lead to a less competitive approximation with respect to compute the full FE solution.
Finally, we highlight that being the problem linear, the coefficients $ \thetaAhNamu{q} $ do not depend on $\rbvecumu $.

\subsubsection{PDE-DNNs for nonaffinely parametrized PDEs}
\begin{figure}[t]
	\centering
	\subfloat[Input location (green) and output location (red)]{\vspace{-1.5cm} \includegraphics[width=0.45\textwidth]{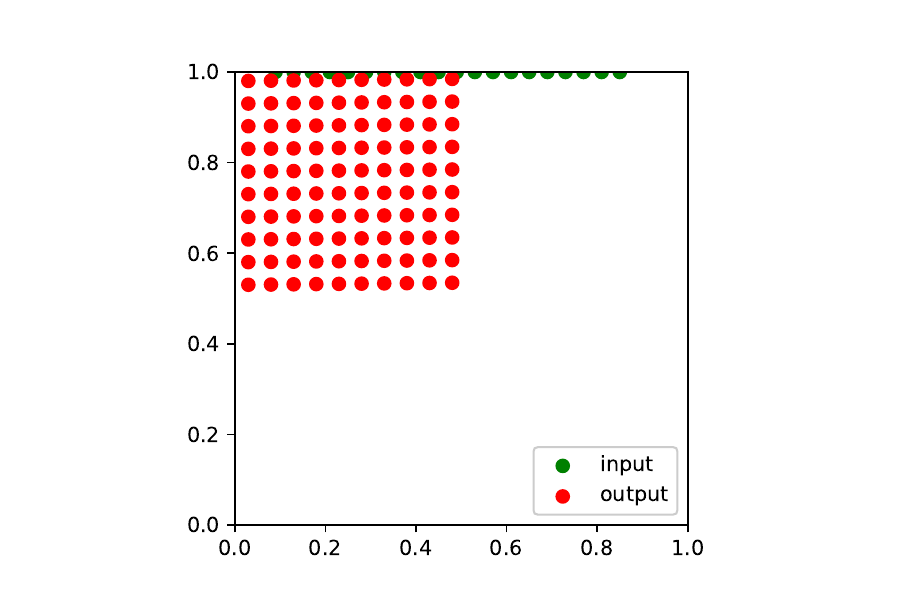}\label{fig:na_points}}
	\subfloat[$\archNaOne $ and $\QaNa = 5$]{\includegraphics[width=0.4\textwidth]{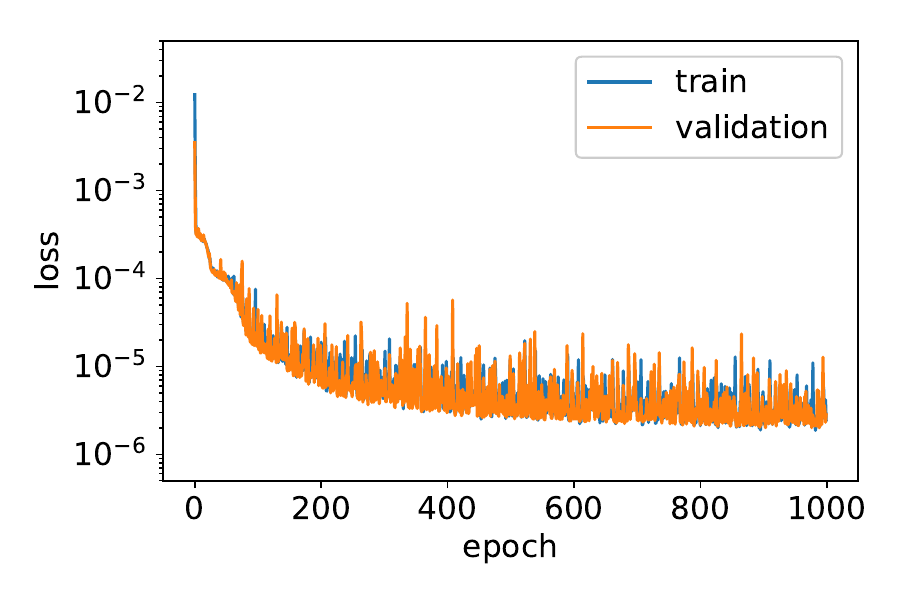} }
	\\
	\subfloat[$\archNaOne $ and $\QaNa = 40$]{\includegraphics[width=0.4\textwidth]{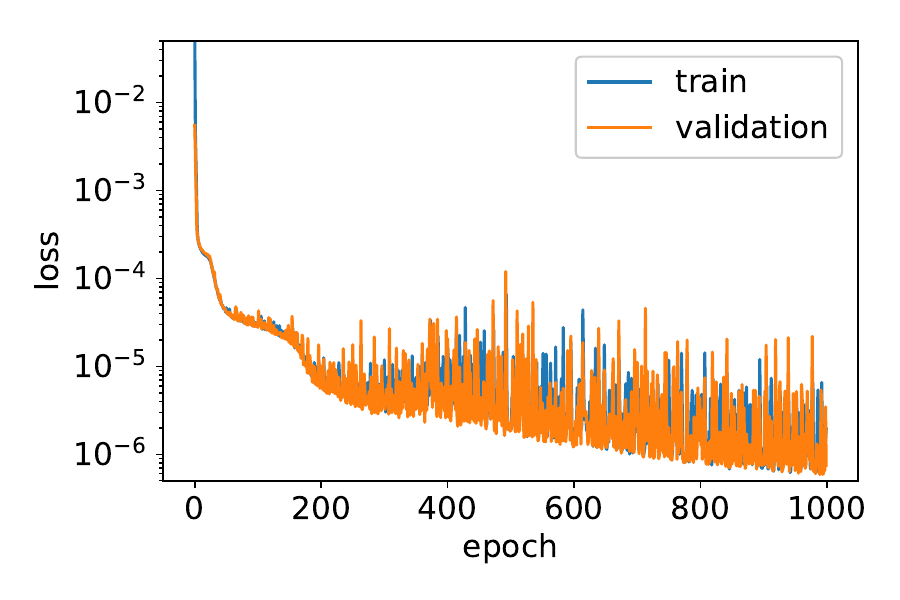} }
	\subfloat[$\archNaThree $ and $\QaNa = 5$]{\includegraphics[width=0.4\textwidth]{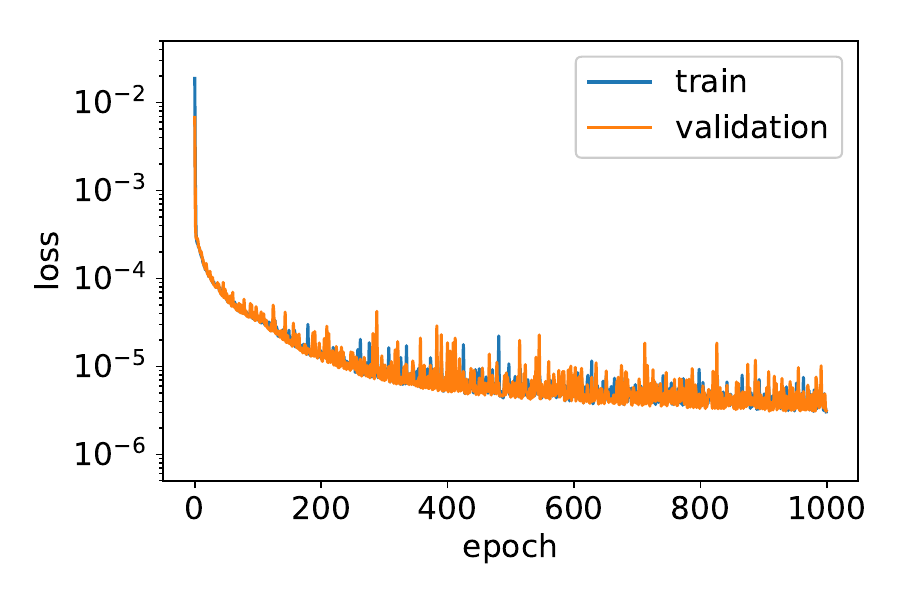} }
	\caption{Location of input and output locations and decay of loss for the training and validation set different architectures and number of affine terms $\QaNa$. }
	\label{fig:trainingNonaffine}
\end{figure}

Given the measurements of the FE solution at $\Nin=40$ physical points located on the boundary $\Gamma_N $ (green in Figure \ref{fig:na_points}), we are interested in predicting the value of the same solution at $\Nout=100$ physical points located in the top left quadrant of the domain (red in Figure \ref{fig:na_points}). To this aim, we employ a PDE-DNN with 4 hidden layers, for which we investigate three different architectures, reported in Table \ref{tab:na_layers}.
\begin{table}
	\centering
	\begin{tabular}{ c | c}
		\toprule
		MLP architecture & layers description\\
		\midrule
		$\archNaOne$ & 4 hidden layers with sizes 1024, 512, 256, 128 \\
		$\archNaTwo$ & 4 hidden layers with all sizes equal to 256 \\
		$\archNaThree$ & 4 hidden layers with all sizes equal to 64 \\
		\bottomrule
	\end{tabular}
	\caption{Architectures employed for Nonaffine test case.}
	\label{tab:na_layers}
\end{table}
The fifth and last layer combines the outputs of the forth layer into the coefficients
$\{\thetaAhNamu{q}\}_{q=1}^{\QaNa} $ and $\{\thetaFNamu{q}\}_{q=1}^{\QfNa} $ of the affine approximations \eqref{eq:affine_approx_na}. These coefficients are used to assemble the RB matrix and right hand side. After solving the small dense linear problem, the layer recovers the values of the
temperature at the output locations.

With respect to the PDE-DNN in Figure \ref{fig:ex1hidden_pde}, $\Dmu{i}$ represent the
following values:
\begin{align*}
\Dmu{i} =
\begin{cases}
\thetaAhNamu{i} \qquad &\operatorname{if} \; i = 1,\dots,\QaNa \\
\thetaFNamu{i-\QaNa} \qquad  &\operatorname{if} \; i = \QaNa+1,\dots,\QaNa+\QfNa.
\end{cases}
\end{align*}
The number of employed affine components is fixed to $\Qf = 10 $ for the right hand side, yielding an accuracy of $\approx 10^{-6} $ for the approximation of the right hand side.
On the other hand, the considered parametrization leads to a large number of terms for the affine approximation of $ \feAhNamu$ to obtain its accurate representation (up to about 40 for an accuracy of $\approx 10^{-5}$), leading to a significant overhead when employing the standard RB method.
In the latter, the values $\{\thetaAhNamu{i}\}_{i=1}^{\QaNa} $ are obtained by solving an interpolation problem, and a large number of affine components is required to obtain an accurate RB solution.
In the PDE-NN framework we are proposing, the values of the functions $\{\thetaAhNamu{i}\}_{i=1}^{\QaNa} $ (and $\{\thetaFNamu{i}\}_{i=1}^{\QfNa} $) are instead optimized by the encoding part of the network, the MLP.
It is important to notice in the case of PDE-NNs the values $\{\thetaAhNamu{i}\}_{i=1}^{\QaNa} $ are not trained to minimize the error between the original FE matrix $\feAhNamu $ and its affine approximation, but rather to minimize the error between the final output of the RB solver and the training output.
We analyze the behavior of our network by varying the number of matrix affine components $\QaNa $, which is set to 1, 2, 3, 4, 5, 10, 20, 40, and used to feed the RB solver.
These numbers of affine components lead to an affine approximation of $\feAhNamu $ that by the standard computation of $\{\thetaAhNamu{i}\}_{i=1}^{\QaNa} $ would lead to a relative error ranging from $1 $ to $10^{-5}$.
We stress that for our PDE-DNNs we employ MDEIM for building an affine basis, but not for computing the corresponding coefficients $\{\thetaAhNamu{i}\}_{i=1}^{\QaNa} $ of the affine approximation.

The training of the PDE-DNN is performed through minimization of the loss in Eq. (2) by Adam algorithm (with learning rate set to $10^{-3}$) for 1000 epochs; the parameters of the optimization method were set after empirical observations.
For training the network, we employ $\Ns = 10000$ samples obtained by computing the FE full solutions of problem \eqref{eq:nonaffineFem} corresponding to randomly sampled values of $\param \in \paramspace $. We use this  set also to construct a RB solver with $\epod = 10^{-4} $. We considered this setting adequate after the analysis carried out in the test case in Section \ref{subsec:affine}; 80\% of the samples are used as training set, whereas the remainder  serves as validation set. Examples of loss decay for the training validation sets corresponding to $\archNaOne $ and $\archNaThree $ architectures and $\QaNa=5 $ and $\QaNa=40 $ affine components are reported in Figure \ref{fig:trainingNonaffine}. As expected, the lower number of trainable parameters in the case of $\archNaThree $ results in a significantly smoother decay of the loss function.

\begin{table}
	\centering
	\begin{tabular}{ c  c  c  c  c  c  c  c  c } 
		\toprule 
		$\QaNa $ & 1 & 2 & 3 & 4 & 5 & 10 & 20 & 40\\ 
		\midrule 
		$\archNaOne $ & 2.19e-02 & 3.59e-02 & 2.55e-02 & 7.31e-03 & 5.15e-03 & 6.16e-03 & 5.58e-03 & 2.59e-03\\ 
		$\archNaTwo $ & 2.22e-02 & 8.58e-03 & 1.18e-02 & 5.46e-03 & 4.87e-03 & 3.11e-03 & 4.55e-03 & 2.36e-03\\ 
		$\archNaThree $ & 2.22e-02 & 2.97e-02 & 6.26e-02 & 7.79e-03 & 5.69e-03 & 4.42e-03 & 4.52e-03 & 3.45e-03\\ 
		RB & 6.84e-01 & 7.63e-01 & 3.43e+00 & 1.34e+00 & 6.52e-01 & 4.98e-02 & 7.27e-03 & 3.32e-04\\ 
		\bottomrule 
	\end{tabular}
	\caption{Errors on the output for the three architectures and the RB method over 1000 test samples randomly selected and different from the ones using in the training phase.}
	\label{tab:na_results}
\end{table}

The trained network is then evaluated on a test set with 1000 samples corresponding to values of the parameters randomly sampled in $\paramspace $ and different from the ones used during the training phase, and we evaluate the error achieved on the value of the solution at the output location by the three proposed architecture as function of the number of affine components predicted by the MLP. We also compare this result by the one computed by a standard RB solver, where the number of affine components is $\QaNa $ and, given the physical parameter $\param$, the coefficients $\{\thetaAhNamu{q}\}_{q=1}^{\QaNa} $ are computed online by solving an interpolation problem in MDEIM.
%The results are reported in Figure \ref{fig:nonaffine_accuracies}, and we compare them to what is obtained by the standalone RB method with the same number of matrix affine components.
As a matter of fact, the RB method requires a fine affine decomposition of $\feAhNamu $ in order to provide a satisfactory result. As observed in the literature, the accuracy of the RB solution strictly depends on the accuracy of the matrix affine approximation: the lower $\QaNa$ is, the less accurate the RB model becomes.
On the contrary, the accuracy of the PDE-DNNs is mildly dependent of the number of affine components $\QaNa $: a minimal amount of $\QaNa = 4 $ affine basis is required, where a plateau is reached.
By employing a larger $\QaNa$ no gain in accuracy is observed.
This different behavior is ascribed to the fact that the PDE-DNN does not necessarily provide the values of $\{\thetaAhNamu{q}\}_{q=1}^{\QaNa} $ that would be computed by the interpolation problem solved by the MDEIM algorithm, and that would provide a poor approximation in the case of small $ \QaNa$. Instead, it computes the coefficients $\{\thetaAhNamu{q}\}_{q=1}^{\QaNa} $ which yield an accurate RB approximation with the given affine basis.
Eventually, this results in the minimization of the loss function.
On average, the PDE-DNN approximation is more accurate in the 100\% of cases when using $\QaNa=1, 2, 3, 4, 5, \,10 $ and in about the 90\% of cases if $\QaNa =20 $. The standalone RB method always results in a better approximation only if $\QaNa=40$. The detailed errors can be found in Table \ref{tab:na_results}.

\begin{figure}
	\centering
	\includegraphics[width=0.5\textwidth]{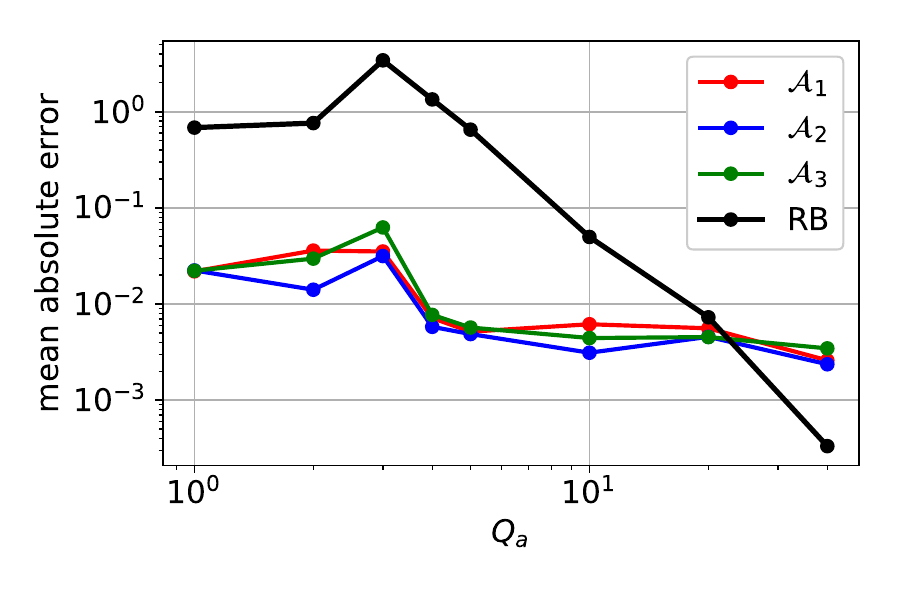}
	\caption{Accuracy obtained by the PDE-DNN with three different architectures of the MLP and the standalone RB method as function the number of matrix affine components $\QaNa $.}
	\label{fig:nonaffine_accuracies}
\end{figure}

\subsection{Steady Navier-Stokes equations}\label{sec:navier_stokes}
\begin{figure}
	\centering
	\subfloat[Regions $ \Omega_s$ (red) and $ \Omega_c$ (green). ]{\includegraphics[width=0.49\textwidth]{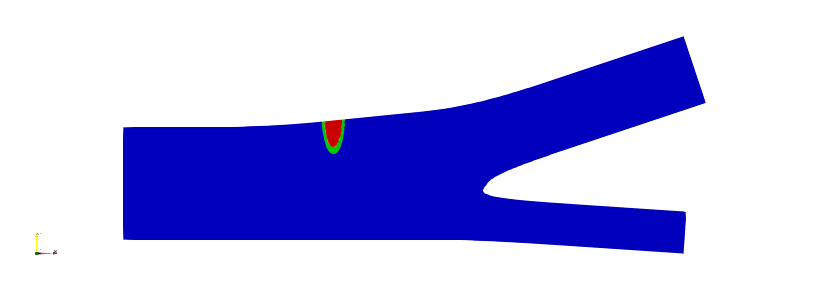}\label{fig:nsMeshesB}}
	\subfloat[FE mesh.]{\includegraphics[width=0.49\textwidth]{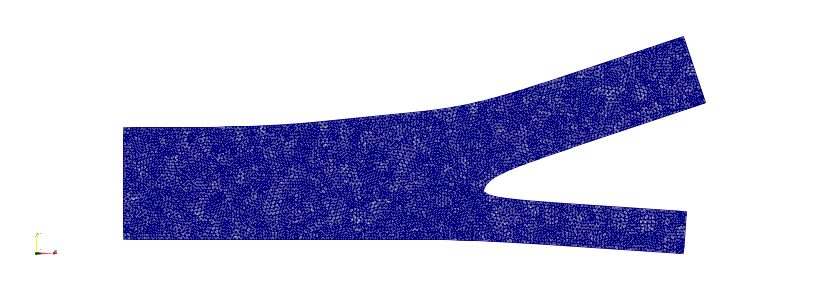}\label{fig:nsMeshesA}}
	\caption{On the left, domain $\Omega$, boundaries, regions $ \Omega_s(r)$, $ \Omega_c(r)$ for $r = 0.2$. On the right, the FE mesh. }
	\label{fig:nsMeshes}
\end{figure}

In this section we introduce the steady Navier-Stokes (NS) equations for a viscous Newtonian incompressible fluid, which model the blood dynamics in artery bifurcations.
Given an open bounded domain $\Omega \subset \real^2 $, shown in Figure \eqref{fig:nsMeshesB}, such that
$\partial\Omega = \Gamma_{\operatorname{in}} \cup \Gamma_{\operatorname{out}} \cup \Gamma_{\operatorname{w}}$ and  $ \mathring{\Gamma}_{w} \cap \mathring{\Gamma}_{in} = \mathring{\Gamma}_{out} \cap \mathring{\Gamma}_{in} =  \mathring{\Gamma}_{w} \cap \mathring{\Gamma}_{out} = \emptyset $ the steady NS equations read as follows:
\begin{align}\label{eq:nsdiff}
\begin{cases}
\displaystyle
\nsu\cdot\nabla\nsu -\nu\Delta\nsu + \nabla\nsp = \vec{0} \qquad & \text{ in } \Omega
\\
\nabla \cdot \nsu = 0	 \qquad & \text{ in } \Omega \\
\mbox{+b.c.} & \text{ on } \partial\Omega
\end{cases}
\end{align}
Here $\nsu = \nsumu  $ and $\nsp = \nspmu $ are the velocity and the pressure fields describing the dynamics of the fluid and
$\nsnu $ denotes the kinematic viscosity, which in our formulation is taken as physical parameter.
We employ a nonhomogeneous Dirichlet inlet condition $\uu=\uu(\vecx)= 4U(0.4-x)^2$, with $U=10 $, on $\Gamma_{\operatorname{in}} $ , homogeneous Dirichlet conditions on the wall $\Gamma_{\operatorname{w}} $ and homogeneoous Neumann conditions on the
outlet boundary $\Gamma_{\operatorname{out}}$.
We define the  Reynolds number $\rey = {L\bar{U}}/{\nsnu}$  as the non-dimensional ratio of convection to diffusion,  where $L$ and $\bar{U} $ are the characteristic length of the domain and velocity of the flow, respectively; here  we deal with laminar flows,  featuring    $\rey \in [1,10^3]$.
We are interested in the velocity dynamics in the presence of stenosis, that is the partial occlusion of the artery before the bifurcation.
In this work, this is modelled by summing an additional reaction term $c(\vecx; r)\nsu $ to the left hand side of the first equation in \eqref{eq:nsdiff}, where
\begin{align*}
c(\vecx; r) =
10^3\mathbb{I}_{\Omega_s(r)} + 10^{-10}\mathbb{I}_{\Omega_c(r)}
\end{align*}
where $\mathbb{I}_A $ is the characteristic function on a set $A\subset \real^2$ and the regions $ \Omega_s(r)$ and $ \Omega_c(r)$ are defined as
\begin{align*}
&\Omega_s = \Omega_s(r) = \{\vecx \in \Omega : 10(\vecx_1-1.5).^2 + (\vecx_2-0.46)^2 <  r^2\}
\\
&\Omega_c = \Omega_c(r) = \{\vecx \in \Omega\backslash\Omega_s(r) : 10(\vecx_1-1.5).^2 + (\vecx_2-0.46)^2 <  r_{\operatorname{max}}^2\},
\end{align*}
where $r_{\operatorname{max}} = 0.25 $ is fixed.
The function $c(\vecx; r)$ defines an ellipsis modelling the presence of an obstacle which represents the stenosis in the region $\Omega_s$. The volume of the stenosis increases by increasing the value of $r$, which is the second parameter for this test case.
A non-zero value is assigned also in the region $\Omega_c $, in order to guarantee that the support of $c(\vecx; r)$ be not parameter dependent.
This ensures that the sparsity pattern of the nonaffinely parametrized FE matrix obtained by discretizing the additional reaction term does not change by changing the value of the parameter $r$. Should this property not hold, we could not apply MDEIM to affinely approximate it.
Eventually, we define the parameter vector as
\begin{align*}
\param = ( \nsnu, r ) \in \paramspace =[ 0.01, 0.1]\times[0, 0.25].
\end{align*}

From \eqref{eq:nsdiff}, we can derive the corresponding variational formulation and subsequently the FE discretization by employing continuous piecewise second order polynomials for the the velocity and continuous piecewise first order polynomials for the pressure. This choice, also known as Taylor-Hood FE basis functions, is motivated by the saddle-point nature of problem  \eqref{eq:nsdiff} and guarantees the well-posedness of the resulting FE nonlinear system.
We refer to, e.g., \cite{elman2005finite,temam1984navier,brezzi1974existence} for further details on the NS equations and their FE discretization. Solving the FE formulation of the NS equations is equivalent to solving a nonlinear system as \eqref{eq:nonlinearSystem} of the following form, where the matrix at the left hand side has a saddle-point structure
\begin{equation}\label{eq:nsBlocks}
%\begin{array}{c}
\begin{bmatrix}
\nsDmu +  \nsCmu  + \nsKmu & \quad \feBh^T \\
\feBh & \mathbf{0}
\end{bmatrix}%, \smallskip
%\end{array}
\begin{bmatrix}
\nsvecu\\
\nsvecp
\end{bmatrix}
=
\begin{bmatrix}
\nsvecfumu\\
\mathbf{0}
\end{bmatrix}.
\end{equation}
Here $\nsvecu \in \real^{\Nhu} $ and $\nsvecp \in \real^{\Nhp}$ are the vector representations of the velocity and pressure FE solutions, $\nsvecfu $ corresponds to the discretization of the nonhomogeneous Dirichlet condition, $\nsDmu$, $\nsCmu,$ $\nsKmu \in \real^{\Nhu\times\Nhu}$ arise from the discretization of the second order differential operator, the nonlinear term and the additional obstruction term in the first equation of \eqref{eq:nsdiff}, respectively, and $\feBh\in \real^{\Nhp\times\Nhu}$ from the discretization of the incompressibility constraint given by the second equation in \eqref{eq:nsdiff}.
We stress that $\feBh $ is not parameter dependent, whereas the nonlinear matrix $\nsCmu $ depends on the parameter by means of the solution $\nsvecu $.
\begin{remark}
	Let us introduce, for $\vec{w}, \vec{u}, \vec{v}$ in a suitable functional space, the trilinear form
	\begin{align}\label{eq:trilinearNs}
	c(\vec{w}, \vec{u}, \vec{v}) = \int_{\Omega} (\vec{w} \cdot \nabla) \vec{u} \cdot \vec{v}d\Omega.
	\end{align}
	The nonlinear matrix $\nsCmu$ arises from the FE discretization of \eqref{eq:trilinearNs} where $\vec{w}=\vec{u}$.
\end{remark}

\begin{figure}
	\centering
	\subfloat[$\param=(0.05, 0.15)$]{\includegraphics[width=0.49\textwidth]{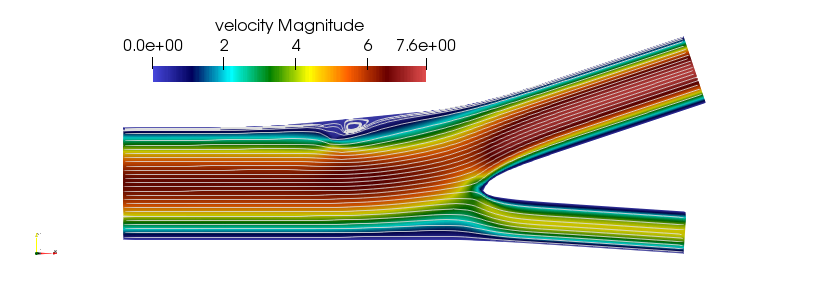}\label{fig:nsSolutionsA}}
	\subfloat[$\param=(0.02, 0.24)$]{\includegraphics[width=0.49\textwidth]{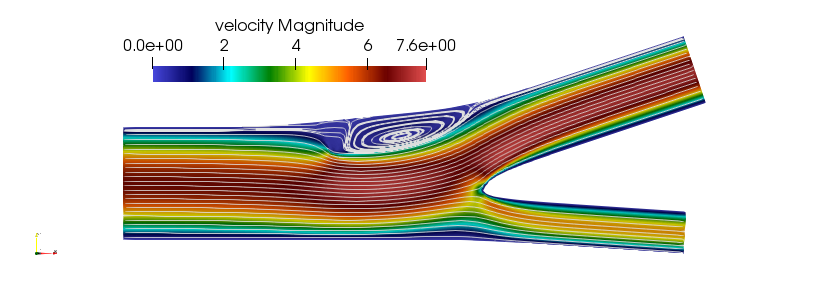}\label{fig:nsSolutionsB}}
	\caption{Velocity solutions of the modified NS problem for $\param=(0.05, 0.15) $ (left) and $\param=(0.02, 0.24)$ (right), where a stenosis is simulated, yielding an obstruction and the subsequent creation of vortices.}
	\label{fig:nsSolutions}
\end{figure}

In the numerical experiment we use a mesh with 15182 triangles, shown in Figure \ref{fig:nsMeshesA}, leading to $\Nhu=59686$ dofs for the velocity and $\Nhp=7592$ for the pressure.
The FE system is then solved with the Newton method (with a tolerance of $10^{-8}$) in about 165.2 seconds, where at each iteration a linear system is solved with a direct method.
An example of regions $\Omega_s $ and $\Omega_c $ is displayed in Figure \ref{fig:nsMeshesB} and the velocity obtained for two choices of the parameter, simulating the absence and the presence of the stenosis, are shown in Figures \ref{fig:nsSolutionsA} and \ref{fig:nsSolutionsB}, respectively.
The larger the value of $r $, the smaller is the volume of central arterial lumen and the larger the corresponding velocity, leading also to the creation of vortices, as it is shown in Figure \ref{fig:nsSolutionsB}.

%\subsubsection{RB methods for the NS equations}
The RB method (and in particular POD) has been largely employed to accelerate the solution of the parametrized steady NS equations, see e.g. \cite{ballarin2015supremizer}. To this aim, as explained in Section \ref{subsec:rbmethod}, we collect velocity solutions of \eqref{eq:nsdiff} for a large enough amount of randomly chosen parameter in $\paramspace$, we perform POD to build the RB projection matrix for velocity $\rbVu\in \real^{\Nhu \times \Nu} $, and approximate $\nsvecu $ in \eqref{eq:nsBlocks}  with $\rbVu \nsrbvecun$, with $\nsrbvecun = \nsrbvecunmu \in \real^{\Nu} $.
Being the columns of $\rbVu $ obtained through POD, i.e. as linear combinations of solutions of \eqref{eq:nsBlocks} for some parameter instances, they are divergence free, that is $\feBh\rbVu = \mathbf{0}$, meaning that the second equation in \eqref{eq:nsBlocks} is automatically satisfied when we substitute $\nsvecu $ with $\rbVu\nsrbvecun$.
Next, we left-multiply the remaining first equation by $\rbVu^T $ to obtain the Galerkin-projected problem. Similarly, $\rbVu^T\feBh^T = \mathbf{0}$, providing us with the following nonlinear RB system where the reduced velocity $\nsrbvecun$ is the only unknown
\begin{align}\label{eq:rb_ns}
%	\rbVu^T\Big( \nsDmu +  \nsCunmu + \nsKmu\Big)\rbVu \nsrbvecun
\Big( \nsrbDmu +  \nsCunmu + \nsrbKmu\Big)\nsrbvecun
=
\nsrbvecfumu.
\end{align}
The definition of the linear reduced operators in \eqref{eq:rb_ns} take advantage of the (approximated) affine decomposition of the corresponding FE arrays, that is,
\begin{align}\label{eq:affine_ns}
\nsrbDmu = \nu\rbVu^T\nsD\rbVu \qquad \nsrbKmu = \sum\limits_{q=1}^{\QkNs} \thetaKhNsmu{q}\rbVu^T\feKhNsq\rbVu,
\qquad \nsvecfumu = \sum\limits_{q=1}^{\QfuNs} \thetaFhNsmu{q}\rbVu^T\feFhNsq,
\end{align}
where $\nsD\in\realhh $ is the FE matrix corresponding to the viscosity-independent Laplace operator, $\{\feKhNsq\}_{q=1}^{\QkNs} \subset \realhh $ and $\{\feFhNsq\}_{q=1}^{\QfuNs} \subset \realh$ are the (M)DEIM basis for the affine approximation of $\nsKmu$ and $\nsvecfumu$, respectively.
For what concerns the nonlinear term $\nsCunmu$, it admits the following (exact) affine decomposition
\begin{align}\label{eq:affine_unj_ns}
\nsCunmu = \sum\limits_{j=1}^{\Nu} \unj \rbVu^T\nsCj\rbVu,
\end{align}
where $\unj$ is the $j $-th component of the RB solution vector, i.e. $\unj = [\nsrbvecun]_j $, and $\nsCj \in \realhh, \, j=1,\dots,\Nu,$ is the FE matrix corresponding to the nonlinear term \eqref{eq:trilinearNs} where $\vec{w} = \vec{\xi}_j $, with $\vec{\xi}_j$ the $j $-th RB function.
The RB problem \eqref{eq:rb_ns} for the steady NS equations is nonlinear and requires to employ a Newton solver and the affine decomposition of the Jacobian and the residual, which are inherited from the ones of the linear and nonlinear operators \eqref{eq:affine_ns} and \eqref{eq:affine_unj_ns}. The need to assemble and solve a nonlinear problem can clearly slow down the efficiency of the RB method.

There are several techniques to obtain a RB problem coupling both the velocity and pressure, which involve the use of the so-called supremizing functions, the interested reader can find more details in \cite{ballarin2015supremizer,rozza2007stability,manzoni2014efficient}.
In this work are not interested to the value of pressure, however, in case of need, it can be reconstructed  from the RB velocity through an additional pressure equation, as in \cite{lovgren2006reduced}.

\subsubsection{PDE-DNNs for parametrized Navier-Stokes PDEs}
We use a setting similar to the one employed in the previous test case and construct a PDE-DNN which has the measurements of the velocity as input and output, with an architecture as \eqref{fig:ex1hidden_pde}, where the RB solver interprets $\Dmu{i}$  as the viscosity $\nu $, the coefficients $\{\thetaKhNsmu{q}\}_{q=1}^{\QkNs} $ and $\{\thetaFhNsmu{q}\}_{q=1}^{\QfuNs} $ of the affine approximations \eqref{eq:affine_ns} and eventually the coefficients $\unj$ of the affine expansion \eqref{eq:affine_unj_ns}.
Since the coefficients $\unj$ have an exponential decay, we decide to truncate them to the first $\QnNs$ coefficients, leading to the approximation
\begin{align}\label{eq:affine_tildeunj_ns}
\nsCunmu \approx \sum\limits_{j=1}^{\QnNs} \unj \rbVu^T\nsCj\rbVu.
\end{align}
To summarize, the last layer processes
\begin{align*}
\Dmu{i} =
\begin{cases}
\nu \qquad &\operatorname{if} \; i = 1\\
\thetaKhNsmu{i-1} \qquad &\operatorname{if} \; i = 2,\dots,\QkNs+1 \\
\thetaFNamu{i-\QkNs-1} \qquad  &\operatorname{if} \; i = \QkNs+2,\dots,\QkNs+\QfuNs+1 \\
\uni{j} \qquad  &\operatorname{if} \; i = \QkNs+\QfuNs+2,\dots,\QkNs+\QfNa+\QnNs, \; j =i-\QkNs-\QfuNs-1 .
\end{cases}
\end{align*}
The choice of predicting the coefficients $\unj$ has evident advantages, since the final layer of our PDE-DNN is called to solve a linear RB problem rather than a non-linear one.
Indeed, we predict the coefficients $\unj$ in the affine decomposition of $\nsCunmu $, instead of considering them as the unknown solution, which is a remarkable advantage when compared to the standard RB method for the NS equations.
As a consequence, also the backpropagation \eqref{eq:backProp} is easily computable.

In the numerical experiments, we have randomly selected a varying number of $\Nin $ input locations and $\Nout $ locations
across the full domain $\Omega $.

\begin{figure}
	\centering
	\subfloat[NS errors at the test locations, $\archNsOne $ architecture. ]{\includegraphics[width=0.45\textwidth]{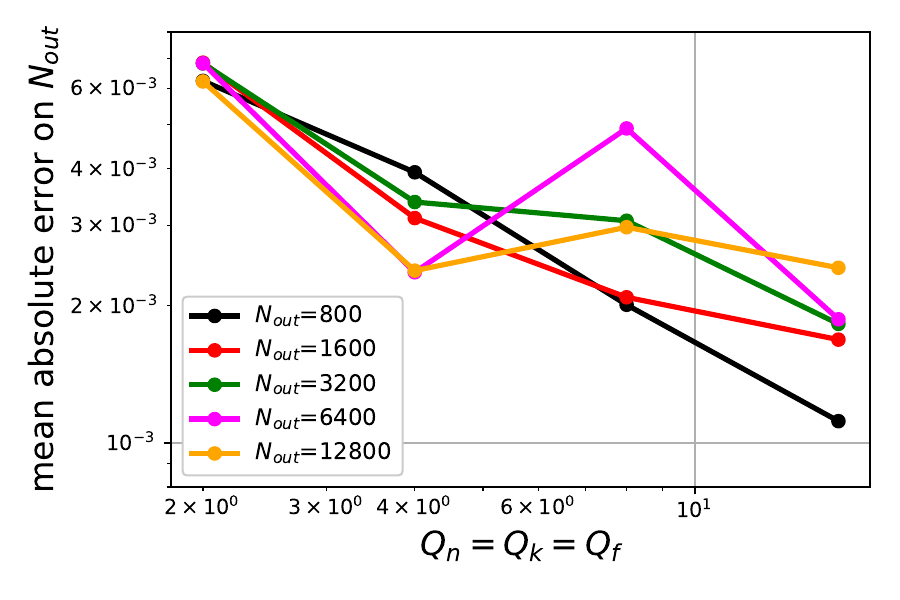}\label{fig:nsErrorsTestSmall}}
	\subfloat[NS errors in $\Omega $, $\archNsOne $ architecture.]{\includegraphics[width=0.45\textwidth]{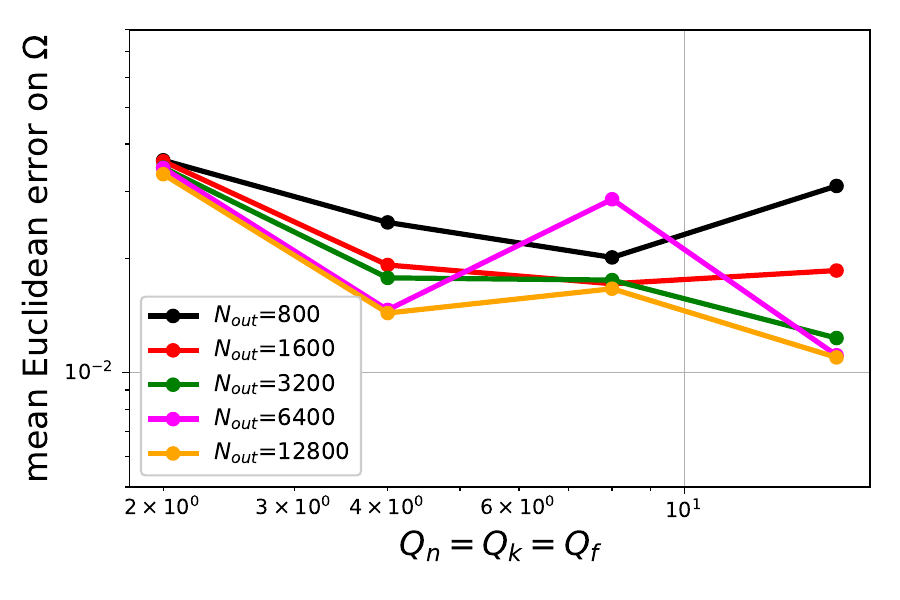}\label{fig:nsErrorsFullSmall}}
	\\
	\subfloat[NS errors at the test locations, $\archNsTwo $ architecture. ]{\includegraphics[width=0.45\textwidth]{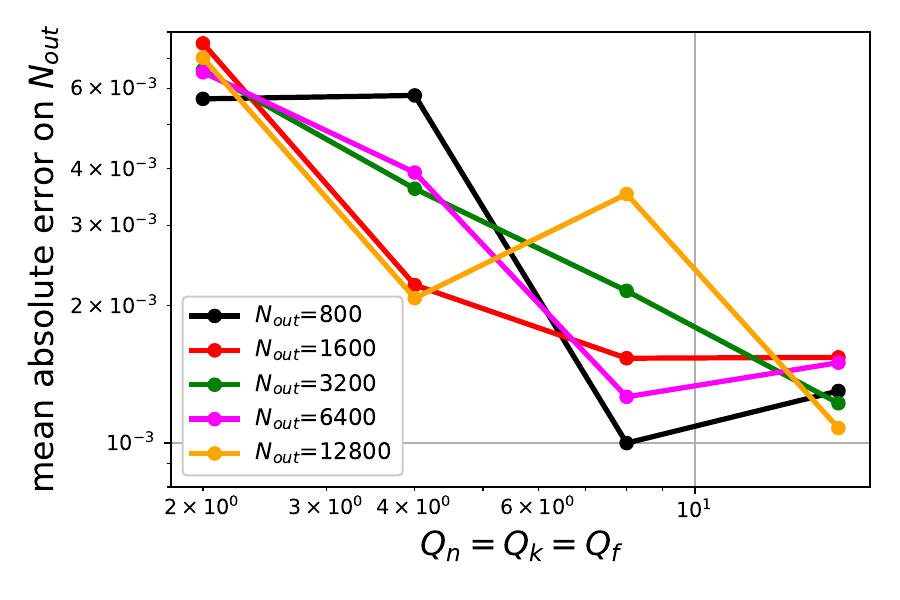}\label{fig:nsErrorsTestNormal}}
	\subfloat[NS errors in $\Omega $, $\archNsTwo $ architecture.]{\includegraphics[width=0.45\textwidth]{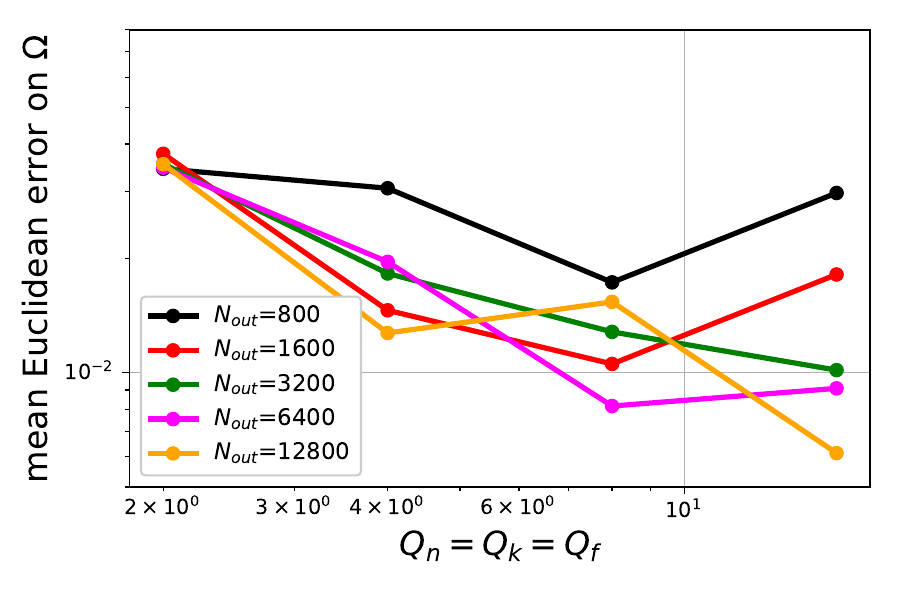}\label{fig:nsErrorsFullNormal}}
	\\
	\subfloat[NS errors at the test locations, $\archNsThree $ architecture. ]{\includegraphics[width=0.45\textwidth]{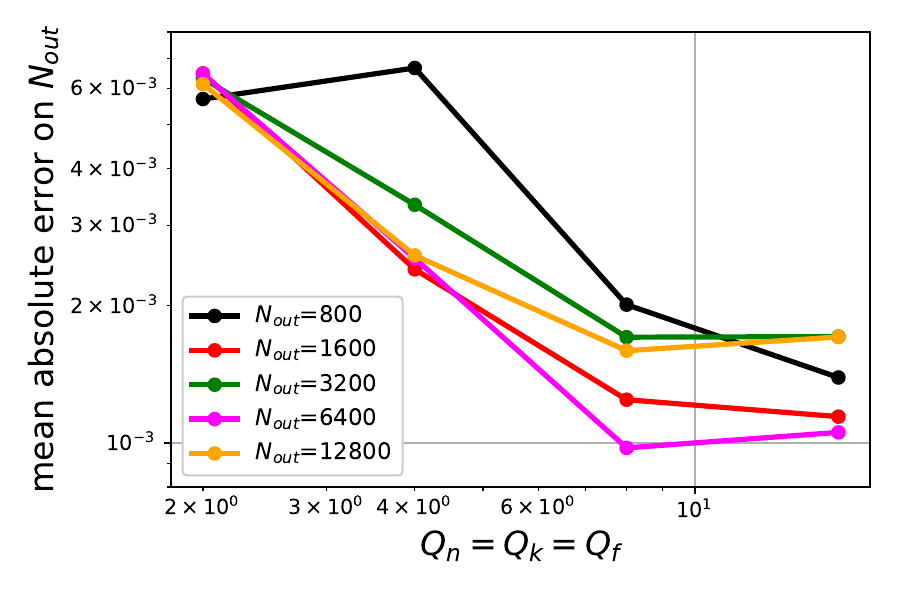}\label{fig:nsErrorsTestLarge}}
	\subfloat[NS errors in $\Omega $, $\archNsThree $ architecture.]{\includegraphics[width=0.45\textwidth]{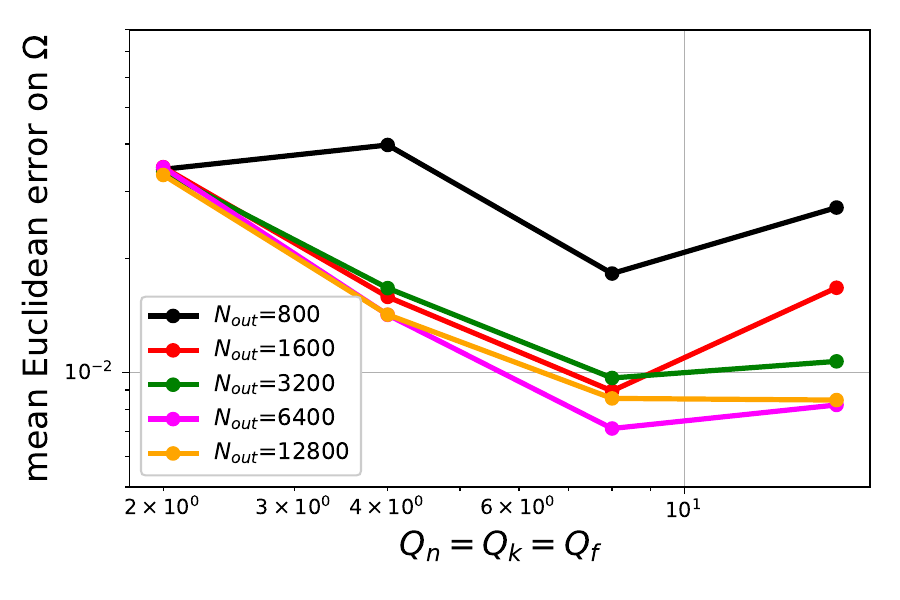}\label{fig:nsErrorsFullLarge}}
	\\
	\subfloat[NS errors at the test locations, $\archNsFour $ architecture. ]{\includegraphics[width=0.45\textwidth]{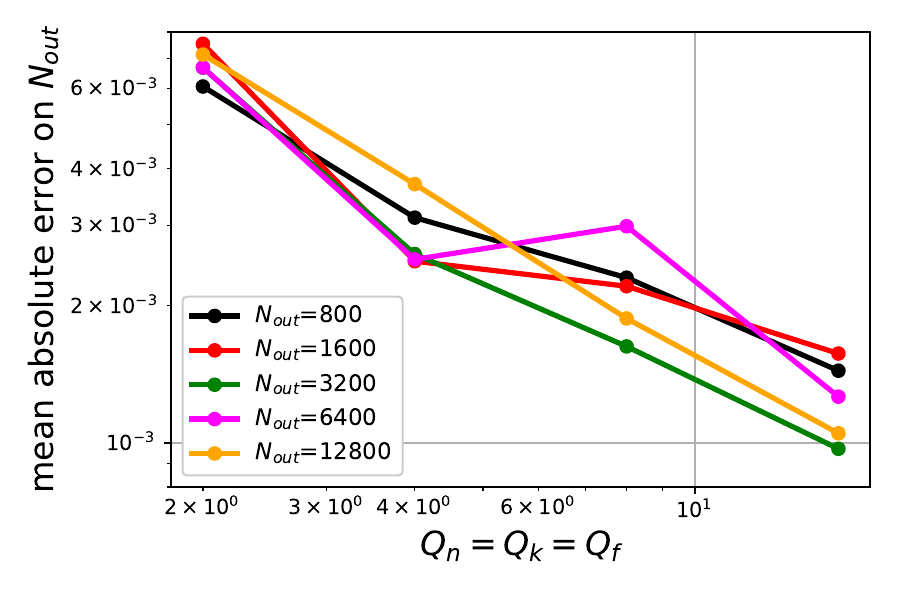}\label{fig:nsErrorsTestLongLarge}}
	\subfloat[NS errors in $\Omega $, $\archNsFour $ architecture.]{\includegraphics[width=0.45\textwidth]{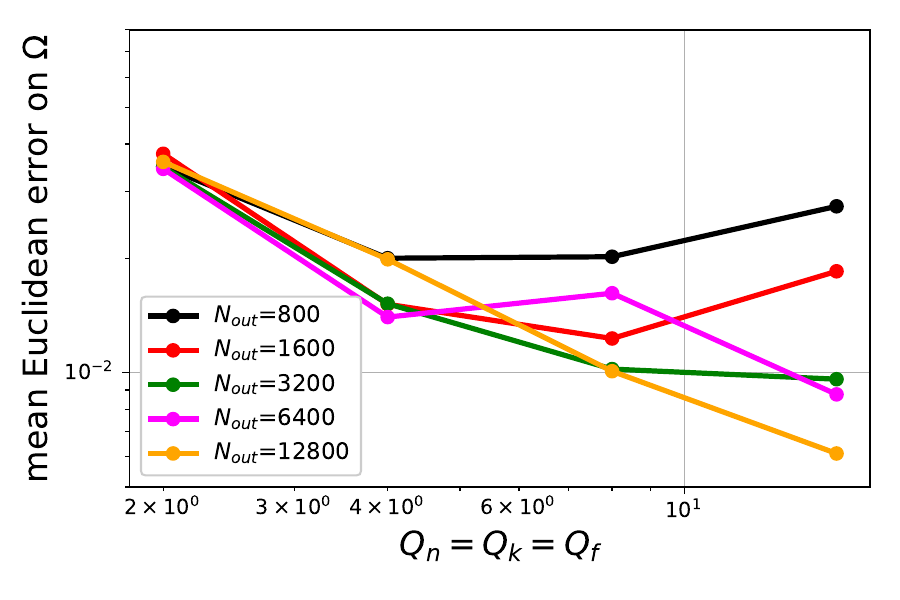}\label{fig:nsErrorsFullLongLarge}}
	
	\caption{On the left, mean error on the test set FE functions at the output locations as function of the number of affine components for the linear (nonaffine) and nonlinear terms. On the right, mean error on the test set FE functions in the full $\Omega $.
		Different lines correspond to different architectures. All results are provided with $\Nin = 100 $.
	}
	\label{fig:nsErrors100}
\end{figure}

We construct the RB problem \eqref{eq:rb_ns} starting from $\ns=19137 $ snapshots for the velocity and with a tolerance $\epod=10^{-4} $, leading to $\Nu=245$ reduced basis functions. Similarly, we employ (M)DEIM with 1000 vector or matrix snapshots to build the affine approximations of $\nsrbvecfumu $ and $\nsrbKmu$; in the experiments we will use $\QkNs=\QfuNs=2, \, 4, \, 8, \, 16 $, which correspond to an affine approximation with a relative accuracy ranging from about 1 to $10^{-2}$.
For the MLP component, we take as architecture the ones reported in Table \ref{tab:ns_layers} and we employ the same $\ns=19137 $ samples obtained by solving the FE problem corresponding to randomly selected instances of $\param$. We train our PDE-DNN with the Adam algorithm and 500 epochs and test against 1000 FE solutions, obtained for randomly selected instances of $\param $, different from the ones employed during the training.
For simplicity, in the following we set the number of predicted affine components for the nonlinear term equal to the one $\QkNs $ for the linear part, thus resulting in $ \QnNs = \QkNs=\QfuNs=2,\, 4,\, 8, \, 16 $.
The errors after training as function of the number of output locations $\Nout $ are reported in Figure \ref{fig:nsErrors100} for the case $\Nin=100$ for many architectures.
\begin{table}
	\centering
	\begin{tabular}{ c | c | c}
		\toprule
		MLP architecture & Number of layers $L$ & Number of nodes per layer \\
		\midrule
		$\archNsOne$   & 4  &  64 \\
		$\archNsTwo$   & 4  &  256 \\
		$\archNsThree$ & 4  &  1024 \\
		$\archNsFour$  & 8  &  1024 \\
		\bottomrule
	\end{tabular}
	\caption{Architectures employed for NS equations.}
	\label{tab:ns_layers}
\end{table}

As a first consideration, we can see that the error has a decreasing trend with the number of affine components used for the nonaffine and nonlinear term, which is however interrupted by few bumps which are motivated by the training with the Adam algorithm with a mini-batch size equal to 128 and the random initialization of the trainable parameters.
Then, we see that employing an architecture with a smaller amount of trainable parameters performs better when also a small number of output locations is employed; as an example we refer to Figure \ref{fig:nsErrorsTestSmall}, where a better result on the test set is achieved with $\Nout = 800 $. 
The behavior can be ascribed to the fact that a smaller number of trainable parameter is not able to capture the complex dynamics considered in this test case, not being able to create a map with too many output points. Yet, an overall better accuracy is achieved when more output points are employed, since in this case the error on a larger portion of the domain is penalized in the loss function used for the training (Figure \ref{fig:nsErrorsFullSmall}).
As one could expect, the better results are provided when the architecture with a larger number of trainable parameter, $\archNsFour $, is taken into account; this situation corresponds to Figure \ref{fig:nsErrorsTestLongLarge}, and achieves the best results in term of accuracy both on the $\Nout $ test output locations and in the full domain. As expected, the larger the number of testing output locations $\Nout $, the better the results on the full domain $\Omega $ (Figure \ref{fig:nsErrorsFullLongLarge}).
Finally, we observe that by using a larger number of output locations, which consists in providing more information to the network during the training phase, yields a more significant use of a larger number of affine operators. Indeed, with $\Nout = 800 $, the accuracy does not significantly depend on the number of affine components; instead, the larger $\Nout$, the better are the improvements given by a larger number of affine components.
The architectures with an intermediate number of trainable parameters, present intermediate results in term of accuracy on both the test set and the full solution.

\begin{figure}
	\centering
	\subfloat[NS errors at the test locations, $\Nin=200 $. ]{\includegraphics[width=0.45\textwidth]{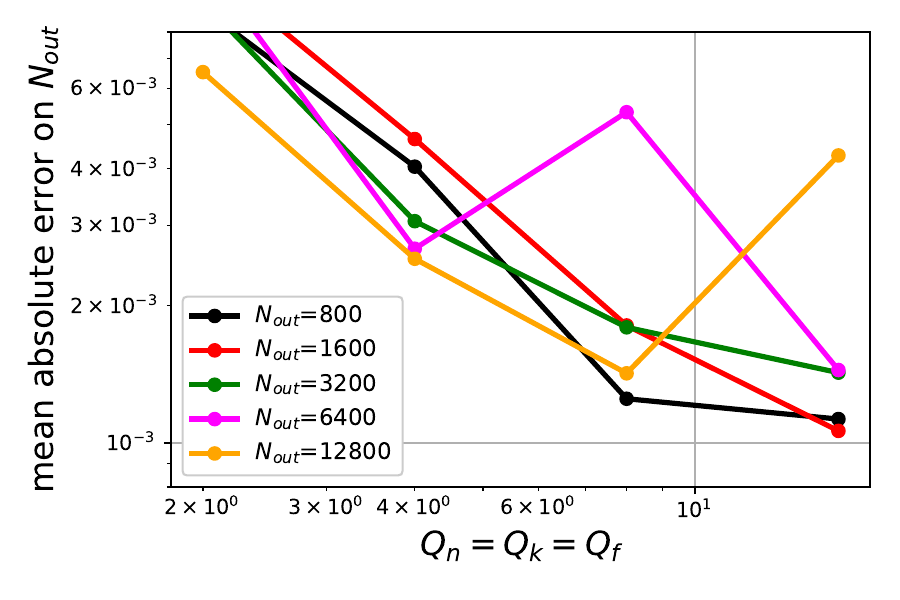}\label{fig:nsErrorsTestLongLarge200}}
	\subfloat[NS errors in $\Omega $, $\Nin=200 $.]{\includegraphics[width=0.45\textwidth]{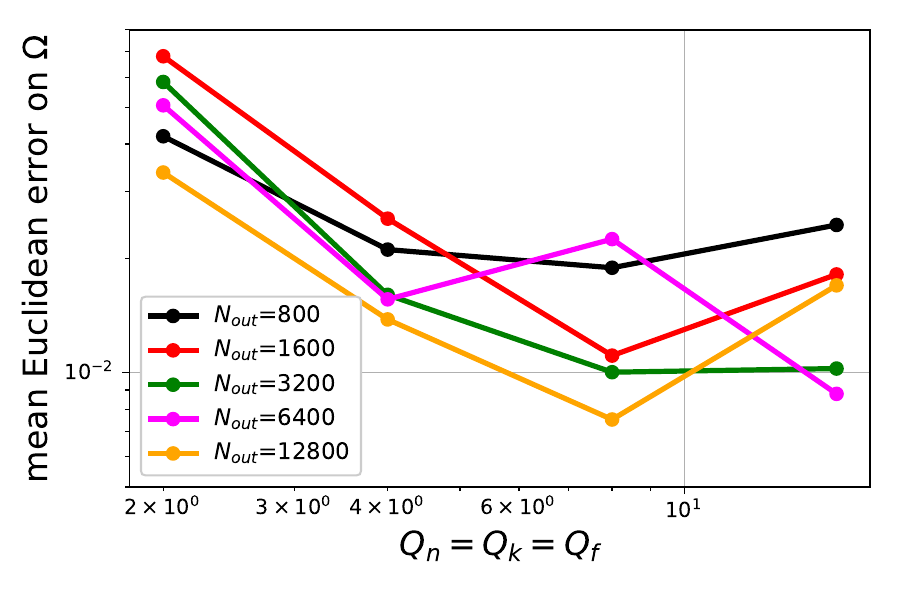}\label{fig:nsErrorsFullLongLarge200}}
	\\
	\subfloat[NS errors at the test locations, $\Nin=400 $. ]{\includegraphics[width=0.45\textwidth]{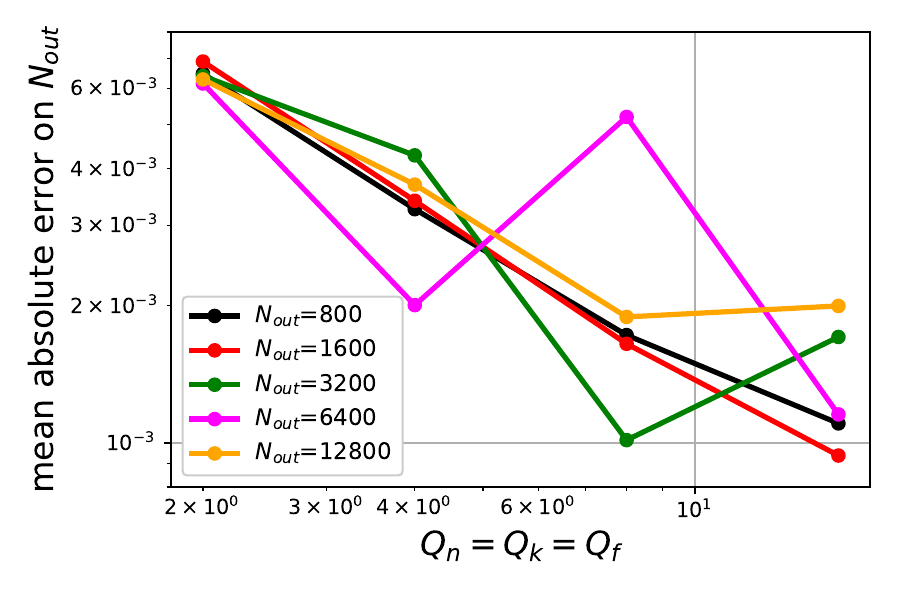}\label{fig:nsErrorsTestLongLarge400}}
	\subfloat[NS errors in $\Omega $, $\Nin=400 $.]{\includegraphics[width=0.45\textwidth]{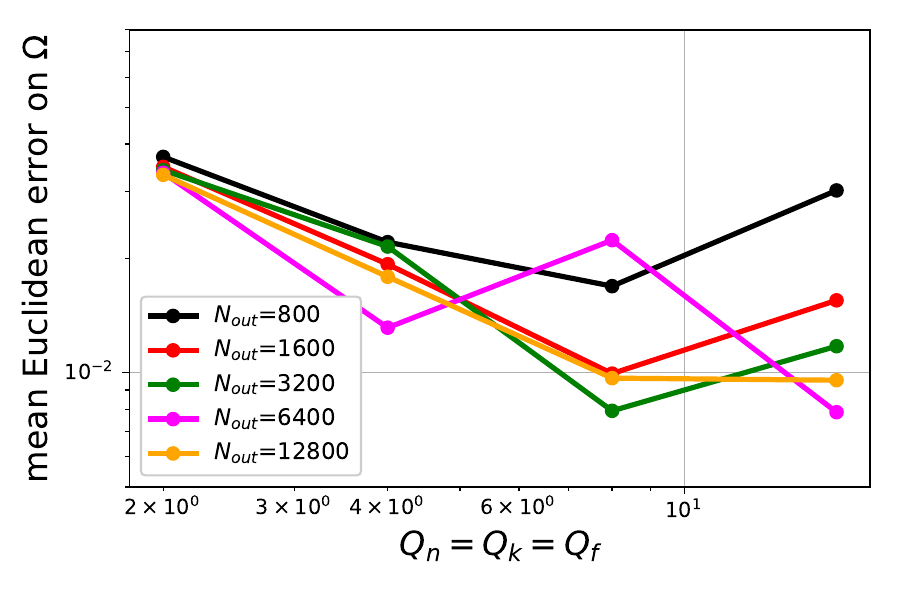}\label{fig:nsErrorsFullLongLarge400}}
	\caption{On the left, mean error on the test set FE functions at the output locations as function of the number of affine components for the linear (nonaffine) and nonlinear terms. On the right, mean error on the test set FE functions in the full $\Omega $.
		The top line refers to $\Nin = 200 $ and the bottom line refers to $\Nin = 400 $. All results are computed with $\archNsFour$ architecture.
	}
	\label{fig:nsErrors200400}
\end{figure}

\begin{figure}
	\centering
	\subfloat[$\Nout = 800 $, $\QnNs=\QkNs=\QfuNs=2 $. ]{\includegraphics[width=0.49\textwidth]{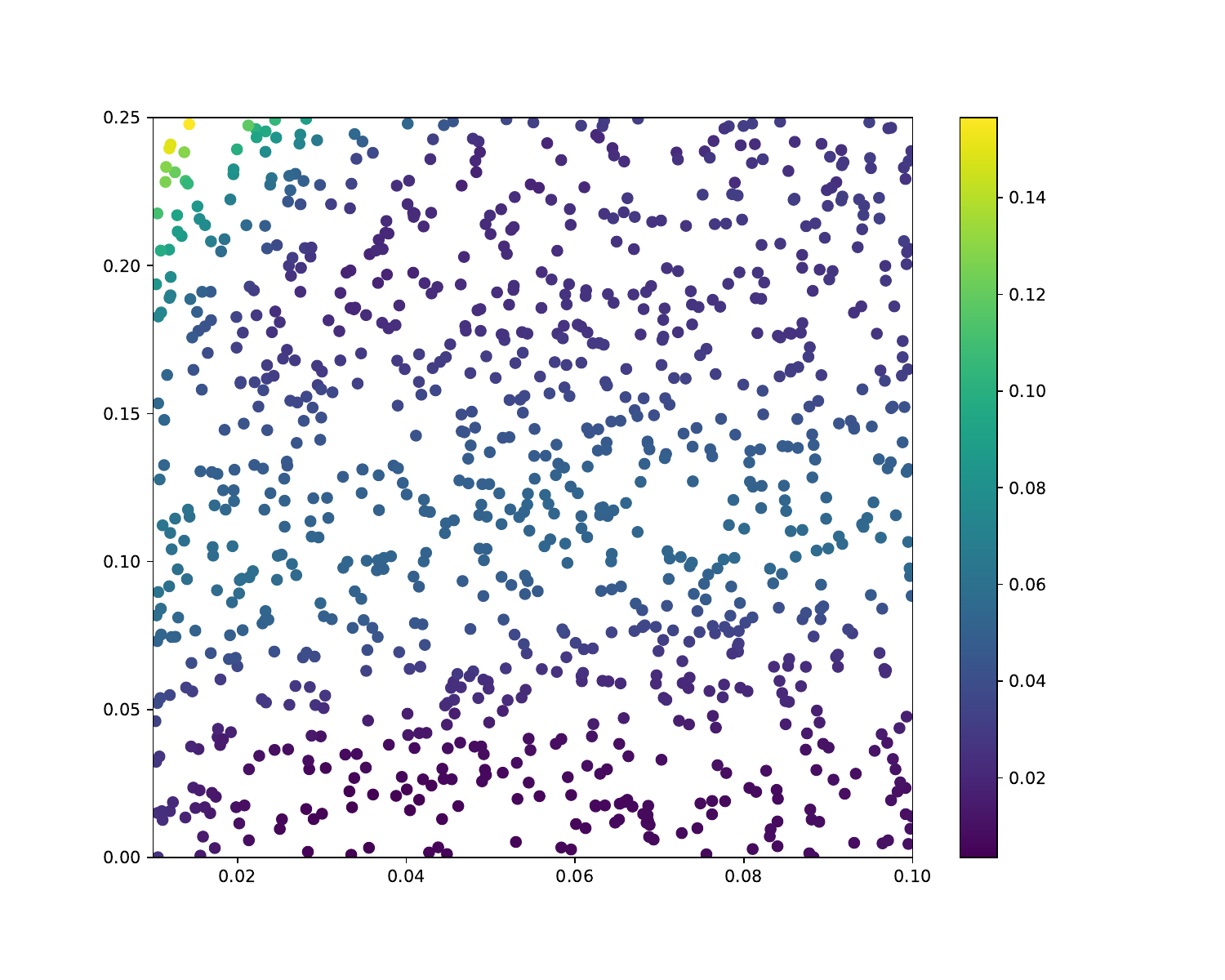}\label{fig:nsScatteredErrors2_800_long}}
	\subfloat[$\Nout = 12800 $, $\QnNs=\QkNs=\QfuNs=2 $.]{\includegraphics[width=0.49\textwidth]{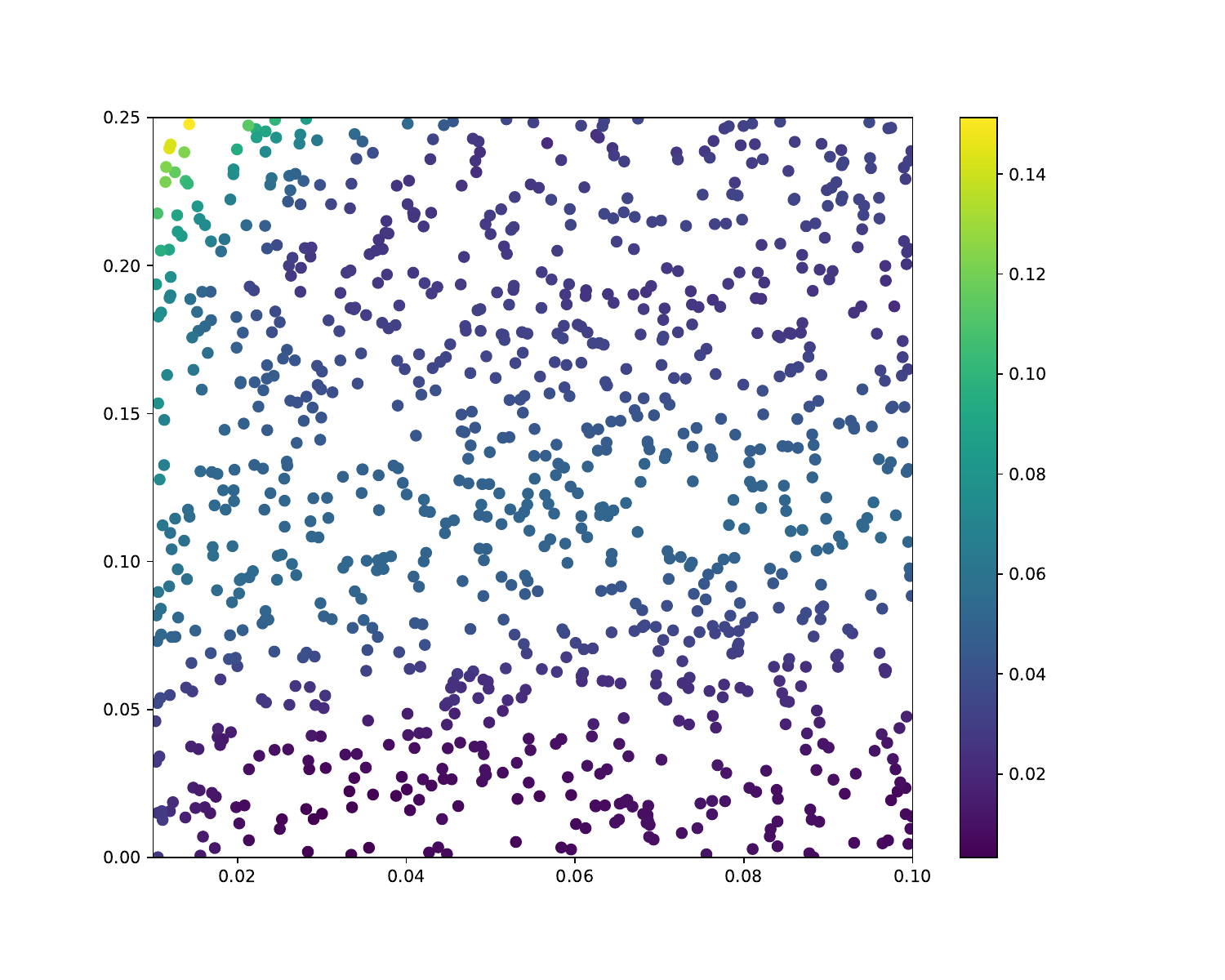}\label{fig:nsScatteredErrors2_12800_long}}
	\\
	\subfloat[$\Nout = 800 $, $\QnNs=\QkNs=\QfuNs=16 $. ]{\includegraphics[width=0.49\textwidth]{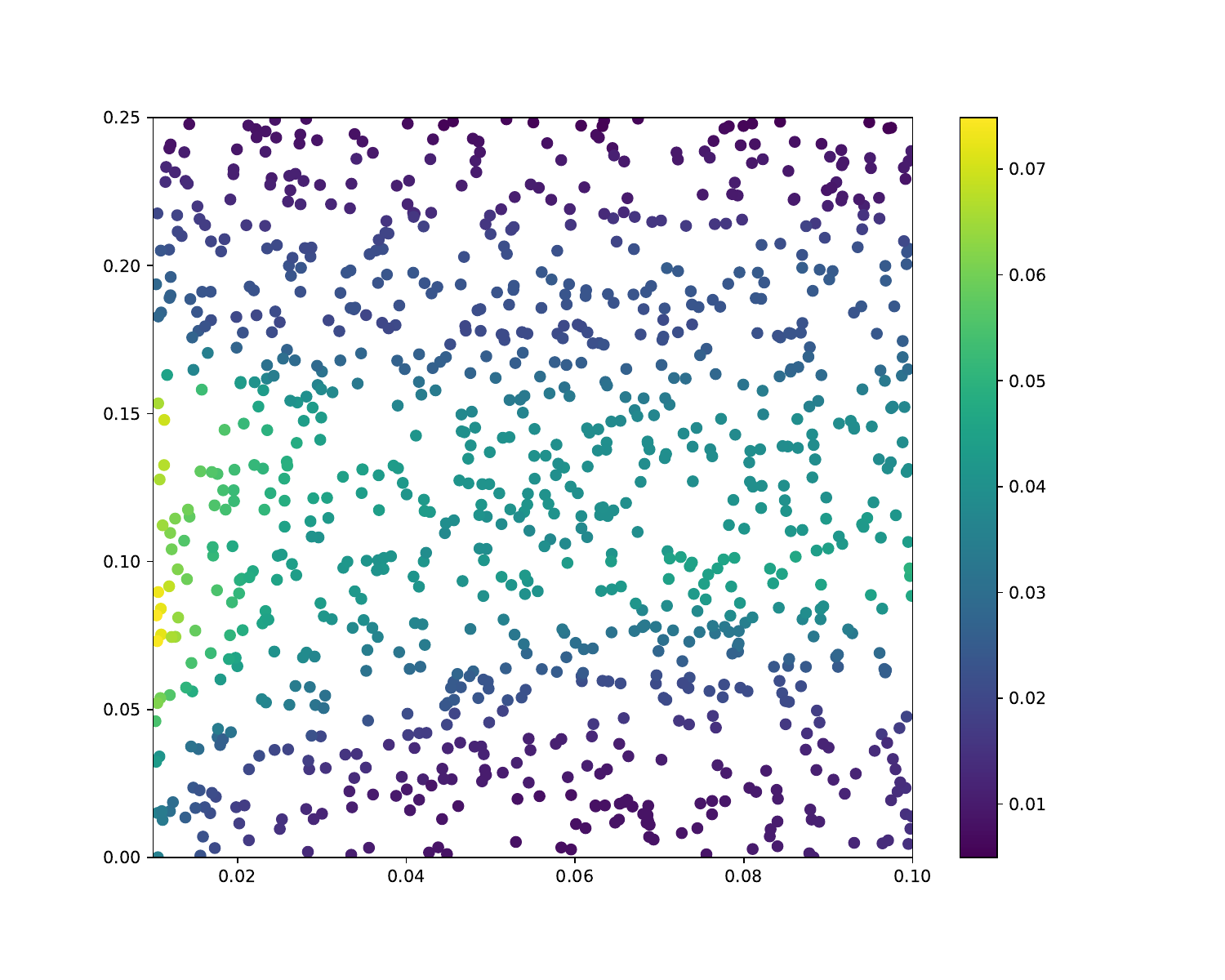}\label{fig:nsScatteredErrors16_800_long}}
	\subfloat[$\Nout = 12800 $, $\QnNs=\QkNs=\QfuNs=16 $.]{\includegraphics[width=0.49\textwidth]{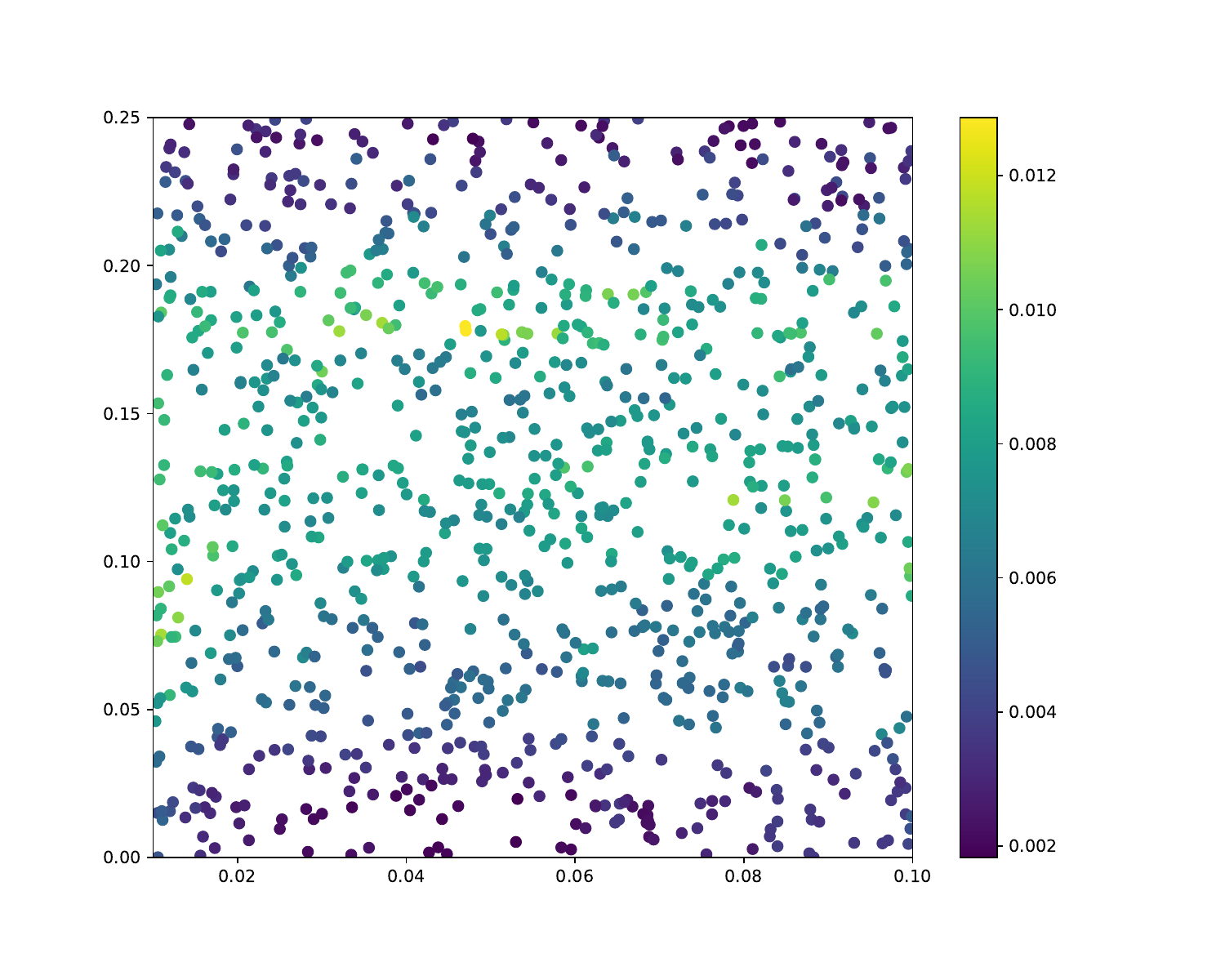}\label{fig:nsScatteredErrors16_12800_long}}
	\caption{Scattered errors (on full $\Omega $) as function of the physical parameters computed with $\archNsFour$ architecture, $\Nin=100 $ and for different configurations of output locations (800 on the left, 12800 on the right) and number of affine components (2 on the top, 16 on the bottom).
	}
	\label{fig:nsScatteredErrors}
\end{figure}

We consider in Figure \ref{fig:nsErrors200400} the cases $\Nin=200, \, 400$, both with the architecture $\archNsFour $; we see that a larger number of input locations does not significantly affect the network in terms of accuracy and we observe the same trends with respect to the number of affine components and the output points as in the test with $\Nin = 100 $.
In Figure \ref{fig:nsScatteredErrors} we report the error for the test parameter values, each point corresponds to a single instance of the parameter and it is colored according to the error committed when reconstructing the solution on the full domain (and compared with the true FE solution). 
The results consider a number of affine component $\QnNs=\QkNs=\QfuNs=2, \, 16 $ (respectively top and bottom line) and a number of output locations $\Nout = 800, \, 12800 $ (respectively left and right).
When using a small number of affine components, i.e.\  $\QnNs=\QkNs=\QfuNs=2$, the error is significantly larger than when using 16 affine components, since the complex dynamics is not captured by such a small number of affine components. 
This is confirmed for example when the fluid is less viscous (small $\nu $) and the occlusion is larger (large $r$), which yield the creation of larger vortices.
Instead, the larger the number of affine components $\QnNs=\QkNs=\QfuNs= 16 $, the smaller the errors and the more effective the use of a larger amount of output locations, up to a maximum error of about 1 \%  across the parameter space.

\section{Conclusions}
\label{sec:conclusions}
In this paper we have proposed a novel way to integrate data and PDE numerical simulation by combining DNNs and RB solvers for the prediction of the solution of a parametrized PDE at some physical points.
The proposed architecture is composed of a MLP followed by a RB solver, where the latter acts as a final nonlinear activation function interpreting the output of the MLP as prediction of parameter dependent quantities - physical parameters (affine case), theta functions of the approximated affine decomposition (nonaffine case) and approximated RB solutions (nonlinear case).
Our architecture recovers an autoencoder structure when the input and output match, in this case the MLP acts as encoder and the RB solver acts as decoder.
We have shown with a wide range of numerical experiments the features of the proposed methodology. Compared to standard DNN, we can obtain as byproduct the solution in the full physical space and, for affine dependencies, the value of the parameter. 
In the nonaffine and nonlinear case, we obtain accurate solutions by employing only a small amount of information compared to what needed by the standard RB method. 
Indeed, for the nonaffine case, only a small number of theta functions is needed by the affine approximation and, in the nonlinear case, we were able to turn the nonlinear problem to a linear one, yet recovering an accurate solution on the full computational domain.

Eventually, the advantages of the methods proposed in this paper can be leveraged for tackling even more complex problem classes. These include, e.g., multiphysics problems, where multiple PDE solvers can be weakly coupled by means of the quantities predicted by the MLP, FE error correctors, to improve the solution given by the considered PDE solver, and, finally, time dependent PDEs, where the (possibly expensive) time marching scheme can take advantage of an extrapolation given by the MLP.

\section*{Acknowledgements}
The authors acknowledge the Swiss National Supercomputing Centre (CSCS) for providing the CPU resources under project ID s796.

\bibliographystyle{abbrv}
\bibliography{biblioSubmission}

\end{document}